# TOROIDAL AND ANNULAR DEHN FILLINGS

Cameron McA. Gordon[1,2] and Ying-Qing Wu[2]

ABSTRACT. Suppose $M$ is a hyperbolic 3-manifold which admits two Dehn fillings $M(r_1)$ and $M(r_2)$ such that $M(r_1)$ contains an essential torus and $M(r_2)$ contains an essential annulus. It is known that $\Delta = \Delta(r_1, r_2) \leq 5$. We will show that if $\Delta = 5$ then $M$ is the Whitehead sister link exterior, and if $\Delta = 4$ then $M$ is the exterior of either the Whitehead link or the 2-bridge link associated to the rational number 3/10. There are infinitely many examples with $\Delta = 3$.

## §1. INTRODUCTION

Let $M$ be a (compact, connected, orientable) 3-manifold with a torus boundary component $T_0$. If $r$ is a *slope* (the isotopy class of an essential simple loop) on $T_0$, then as usual we denote by $M(r)$ the 3-manifold obtained by $r$-*Dehn filling* on $M$, that is, attaching a solid torus $J$ to $M$ along $T_0$ in such a way that $r$ bounds a disk in $J$.

We shall say that a 3-manifold $M$ is *hyperbolic* if $M$ with its boundary tori removed has a complete hyperbolic structure of finite volume with totally geodesic boundary. If $M$ has non-empty boundary, then Thurston's Geometrization Theorem for Haken manifolds [Th1] [Th2] asserts that $M$ is hyperbolic if and only if it contains no essential surface of non-negative euler characteristic, i.e., sphere, disk, torus or annulus.

There has been a considerable amount of investigation of the question of when a Dehn filling on a hyperbolic 3-manifold fails to be hyperbolic; if $\partial M - T_0$ is non-empty then, as noted above, this is equivalent to the existence of an essential surface of non-negative euler characteristic in the filled manifold. So suppose that $M$ is hyperbolic, and that $M(r_1)$ and $M(r_2)$ contain essential surfaces $\widehat{F}_1$ and $\widehat{F}_2$ respectively, each of which is a sphere, disk, torus or annulus. Then for each of the ten possible pairs of surfaces, one can find upper bounds on $\Delta(r_1, r_2)$, the minimal geometric intersection number of $r_1$ and $r_2$. In all cases, the best possible bounds are now known; see [GW] for a fuller discussion. It is also of interest to examine the examples $(M; r_1, r_2)$ that realize values of $\Delta(r_1, r_2)$ at or near the maximal value; it seems to be the case that these are usually quite special.

1991 *Mathematics Subject Classification*. Primary 57N10.
[1] Partially supported by NSF grant #DMS 9626550.
[2] Research at MSRI supported in part by NSF grant #DMS 9022140.

Typeset by $\mathcal{AMS}$-TEX

In the present paper we consider the case where $\widehat{F}_1$ is an annulus and $\widehat{F}_2$ is a torus. Here it is known that $\Delta(r_1, r_2) \leq 5$ [Go, Theorem 1.3]. To state our main result, let $L_1, L_2, L_3$ be the links shown in Figure 7.1 (a), (b) and (c), respectively. These are the Whitehead link, the 2-bridge link associated to the rational number $3/10$, and the Whitehead sister link. Let $M_i$ be the exterior of $L_i$, $i = 1, 2, 3$. Also, let $J_\alpha$ denote the attached solid torus in $M(r_\alpha)$, $\alpha = 1, 2$.

**Theorem 1.1.** *Suppose $M$ is a hyperbolic 3-manifold, such that $M(r_1)$ contains an essential annulus, $M(r_2)$ contains an essential torus, and $\Delta(r_1, r_2) \geq 4$. Then either*
  *(1) $\Delta(r_1, r_2) = 4$ and $M$ is homeomorphic to $M_1$ or $M_2$; or*
  *(2) $\Delta(r_1, r_2) = 5$ and $M$ is homeomorphic to $M_3$.*
  *Moreover, $M_i$ is hyperbolic, $i = 1, 2, 3$, and admits two Dehn fillings $M_i(r_1)$ and $M(r_2)$, such that*
  *(i) $\Delta(r_1, r_2) = 4$ for $i = 1, 2$, and $\Delta(r_1, r_2) = 5$ for $i = 3$;*
  *(ii) $M_1(r_1)$ contains an essential annulus and an essential torus, each intersecting $J_1$ once. All the other $M_i(r_\alpha)$ contain an essential annulus and essential torus, each intersecting $J_\alpha$ twice.*

This follows immediately from Theorems 6.10(1) and 7.5.

In particular, the upper bound of 5 on $\Delta(r_1, r_2)$ is realized by the exterior of the Whitehead sister link. The existence of annular and toroidal fillings on this manifold with $\Delta(r_1, r_2) = 5$ has also been proved by Miyazaki and Motegi [MM].

We also show the following.

**Theorem 7.8.** *There are infinitely many hyperbolic manifolds $M_s$, which admit two Dehn fillings $M_s(r_1), M_s(r_2)$, each containing an essential annulus and an essential torus, with $\Delta(r_1, r_2) = 3$.*

We now sketch the proof and describe the organization of the paper. Throughout, $\alpha, \beta$ will denote the numbers 1 or 2, with the convention that if they both appear, then $\{\alpha, \beta\} = \{1, 2\}$.

We may assume that the surface $\widehat{F}_\alpha$ meets the solid torus $J_\alpha$ in $n_\alpha$ meridional disks, and that $n_\alpha$ is minimal over all choices of $\widehat{F}_\alpha$. The arcs of intersection of the corresponding punctured annulus and punctured torus in $M$ define labeled graphs $\Gamma_\alpha \subset \widehat{F}_\alpha$, with $n_\alpha$ vertices, $\alpha = 1, 2$. We assume that $\Delta(r_1, r_2) = 4$ or 5. The proof consists of a detailed analysis of such graphs $\Gamma_\alpha$, and the elimination of all but four possible pairs, which are then shown to correspond to the examples listed in Theorem 1.1.

In Section 2 we give some definitions, and establish some properties of the intersection graphs $\Gamma_\alpha$, which will be used throughout the paper.

Section 3 deals with the case that one of the graphs has a single vertex. It is shown that if $n_\alpha = 1$ then $\Delta(r_1, r_2) = 4$, $n_\beta = 2$, and $M$ is the exterior of the Whitehead link.

Sections 4 and 5 deal with the case that one of the graphs has two vertices. In Section 4 it is shown that if $n_2 = 2$ then the graph $\Gamma_1$ is 2-separable, i.e., there is an essential sub-annulus $A'$ of $\widehat{F}_1$, which contains two vertices of $\Gamma_1$ and has boundary disjoint from $\Gamma_1$. Note that this is also (trivially) true if $n_1 = 2$. In Section 5 it is shown that if $\Gamma_1$ is 2-separable then $M$ must be one of the three manifolds $M_i$ in Theorem 1.1.

Section 6 deals with the generic case. It is shown that if both graphs have at least three vertices, then $\Delta(r_1, r_2) \leq 3$. This completes the proof of the first part of Theorem 1.1. We then prove in Section 7 that the manifolds $M_i$ in Theorem 1.1 do admit two annular and toroidal Dehn fillings with $\Delta(r_1, r_2) = 4$ or $5$, and that there are infinitely many manifolds which admit two annular and toroidal Dehn fillings with $\Delta(r_1, r_2) = 3$.

## §2. Preliminary Lemmas

Throughout this paper, we will fix a hyperbolic 3-manifold $M$, with a torus $T_0$ as a boundary component. Let $r_1$, $r_2$ be slopes on $T_0$. Denote by $\Delta = \Delta(r_1, r_2)$ the minimal geometric intersection number between $r_1$ and $r_2$. Recall that $\alpha, \beta$ denote the numbers 1 or 2, with the convention that if they both appear, then $\{\alpha, \beta\} = \{1, 2\}$.

**Assumption.** We will always assume that $\widehat{F}_1$ is an essential annulus in $M(r_1)$, and $\widehat{F}_2$ is an essential torus in $M(r_2)$. Moreover, we assume that $M(r_\alpha)$ contains no essential sphere, i.e. it is irreducible.

The last assumption is necessary in our proof, as some of the properties stated in Lemmas 2.3–2.5 are not true without it. However, it is not necessary in the hypothesis of Theorem 6.10, because if $M(r_1)$ is reducible, then by [Wu1] or [Oh] we have $\Delta \leq 3$. Similarly, if $M(r_2)$ is reducible it follows from [Wu3] that $\Delta \leq 2$.

Denote by $F_\alpha$ the punctured surface $\widehat{F}_\alpha \cap M_\alpha$. Let $n_\alpha$, $\alpha = 1, 2$, be the number of boundary components of $F_\alpha$ on $T_0$. Choose $\widehat{F}_\alpha$ in $M(r_\alpha)$ so that $n_\alpha$ is minimal among all essential surfaces homeomorphic to $\widehat{F}_\alpha$. Minimizing the number of components of $F_1 \cap F_2$ by an isotopy, we may assume that $F_1 \cap F_2$ consists of arcs and circles which are essential on both $F_\alpha$. Denote by $J_\alpha$ the attached solid torus in $M(r_\alpha)$. Let $u_1, \ldots, u_{n_1}$ be the disks that are the components of $\widehat{F}_1 \cap J_1$, labeled successively when traveling along $J_1$. Similarly let $v_1, \ldots, v_{n_2}$ be the disks in $\widehat{F}_2 \cap J_2$. Let $\Gamma_\alpha$ be the graph on $\widehat{F}_\alpha$ with the $u_i$'s or $v_j$'s, respectively, as (fat) vertices, and the arc components of $F_1 \cap F_2$ as edges. The minimality of the number of components in $F_1 \cap F_2$ and the minimality of $n_\alpha$ imply that $\Gamma_\alpha$ has no trivial loops, and that each disk face of $\Gamma_\alpha$ in $\widehat{F}_\alpha$ has interior disjoint from $F_\beta$.

If $e$ is an edge of $\Gamma_1$ with an endpoint $x$ on a fat vertex $u_i$, then $x$ is labeled $j$ if $x$ is in $u_i \cap v_j$. In this case $e$ is called a *j-edge at* $u_i$, and in $\Gamma_2$ it is an $i$-edge at $v_j$. The labels in $\Gamma_1$ are considered as integers mod $n_2$. In particular, $n_2 + 1 = 1$. When going around $\partial u_i$, the labels of the endpoints of edges appear as $1, 2, \ldots, n_2$ repeated $\Delta$ times. Label the endpoints of edges in $\Gamma_2$ similarly.

Each vertex of $\Gamma_\alpha$ is given a sign according to whether $J_\alpha$ passes $\widehat{F}_\alpha$ from the positive side or negative side at this vertex. Two vertices of $\Gamma_\alpha$ are *parallel* if they have the same sign, otherwise they are *antiparallel*. Note that if $\widehat{F}_\alpha$ is a separating surface, then $n_\alpha$ is even, and $v_i, v_j$ are parallel if and only if $i, j$ have the same parity.

We use $val(v, G)$ to denote the valency of a vertex $v$ in a graph $G$. If $G$ is clear from the context, we denote it by $val(v)$.

**Definition 2.1.** An edge of $\Gamma_\alpha$ is a *positive edge* if it connects parallel vertices. Otherwise it is a *negative edge*.

When considering each family of parallel edges of $\Gamma_\alpha$ as a single edge $\widehat{e}$, we get the *reduced graph* $\widehat{\Gamma}_\alpha$ on $\widehat{F}_\alpha$. It has the same vertices as $\Gamma_\alpha$.

We use $\Gamma_\alpha^+$ (resp. $\Gamma_\alpha^-$) to denote the subgraph of $\Gamma_\alpha$ whose vertices are the vertices of $\Gamma_\alpha$ and whose edges are the positive (resp. negative) edges of $\Gamma_\alpha$. Similarly for $\widehat{\Gamma}_\alpha^+$ and $\widehat{\Gamma}_\alpha^-$.

A cycle in $\Gamma_\alpha$ consisting of positive edges is a *Scharlemann cycle* if it bounds a disk with interior disjoint from the graph, and all the edges in the cycle have the same pair of labels $\{i, i+1\}$ at their two endpoints, called the *label pair* of the Scharlemann cycle. A pair of edges $\{e_1, e_2\}$ is an *extended Scharlemann cycle* if there is a Scharlemann cycle $\{e_1', e_2'\}$ such that $e_i$ is parallel and adjacent to $e_i'$. There is an obvious extension of this definition to extended Scharlemann cycles of arbitrary length, but this is enough for our purposes.

A subgraph $G$ of a graph $\Gamma$ on a surface $F$ is *essential* if it is not contained in a disk in $F$. The following lemma contains some common properties of the graphs $\Gamma_\alpha$.

**Lemma 2.2.** (Properties of $\Gamma_\alpha$)

*(1)* (The Parity Rule) *An edge $e$ is a positive edge in $\Gamma_1$ if and only if it is a negative edge in $\Gamma_2$.*

*(2) A pair of edges cannot be parallel on both $\Gamma_1$ and $\Gamma_2$.*

*(3) If $\Gamma_\alpha$ has a set of $n_\beta$ parallel negative edges, then on $\Gamma_\beta$ they form mutually disjoint essential cycles of equal length.*

*(4) If $\Gamma_\alpha$ has a Scharlemann cycle, then $\widehat{F}_\beta$ is separating, and $n_\beta$ is even.*

*(5) If $\Gamma_\alpha$ has a Scharlemann cycle $\{e_1, \ldots, e_k\}$, then on $\Gamma_\beta$ the subgraph consisting of $e_1, \ldots, e_k$ and their end vertices, is essential.*

*(6) If $n_\beta > 2$, then $\Gamma_\alpha$ contains no extended Scharlemann cycle.*

*Proof.* (1) can be found in [CGLS, Page 279].

(2) is [Go, Lemma 2.1]. It was shown that if a pair of edges is parallel on both graphs, then the manifold $M$ contains a $(1, 2)$ cable space, which is impossible because $M$ is assumed hyperbolic. See also [GLi].

(3) is [Go, Lemma 2.3] and [GLi, Proposition 1.3]. A set of $n_\beta$ parallel edges with ends $v, v'$ determine a permutation $\varphi$ on the set of labels $\{1, \ldots, n_\beta\}$ such that an edge with label $i$ at $v$ has label $\varphi(i)$ at $v'$. The edges a collection of mutually disjoint cycles of length $L$ in $\Gamma_\beta$, where $L$ is the length of the orbits of $\varphi$. It was shown that if some cycle is not essential then $M$ is cabled, hence non-hyperbolic.

(4) and (5) follow from the proof of [CGLS, Lemma 2.5.2]. It was shown that using the disk bounded by the Scharlemann cycle one can find another surface $\widehat{F}_\beta'$ in $M(r_\beta)$ which has fewer intersections with $J_\beta$, and is cobordant to $\widehat{F}_\beta$, so if $\widehat{F}_\beta$ were nonseparating then $\widehat{F}_\beta'$ would still be essential, which would contradict the minimality of $n_\beta$. If the subgraph $G$ consisting of $e_1, \ldots, e_k$ and their end vertices is inessential then one can show that $M(r_\beta)$ would be reducible, which contradicts our assumption.

(6) If $\widehat{F}_\beta$ is an annulus or torus, this is [Wu3, Lemma 5.4(3)]. If $\widehat{F}_\beta$ is a torus, see [BZ, Lemma 2.9] or [GLu2, Theorem 3.2]. □

Let $\widehat{e}$ be a collection of at least $n_\alpha$ parallel negative edges on $\Gamma_\beta$, connecting $v_1$ to $v_2$. Then $\widehat{e}$ defines a permutation $\varphi : \{1, \ldots, n_\alpha\} \to \{1, \ldots, n_\alpha\}$, such that an edge $e$ in $\widehat{e}$ has label $k$ at $v_1$ if and only if it has label $\varphi(k)$ at $v_2$. Call $\varphi$ the *permutation associated to $\widehat{e}$*. If we interchange $v_1$ and $v_2$, then the permutation is $\varphi^{-1}$; hence it is well defined up to inversion.

**Lemma 2.3.** (Properties of $\Gamma_1$) *(1) If a family of parallel edges in $\Gamma_1$ contains more than $n_2$ edges, and if $n_2 > 2$, then all the vertices of $\Gamma_2$ are parallel, and the permutation associated to this family is transitive.*

*(2) If $\Gamma_1$ contains two Scharlemann cycles with disjoint label pairs $\{i, i+1\}$ and $\{j, j+1\}$, then $i \equiv j \mod 2$.*

*(3) Suppose $n_2 > 2$. Then a family of parallel positive edges in $\Gamma_1$ contains at most $n_2/2 + 2$ edges, and if it does contain $n_2/2 + 2$ edges, then $n_2 \equiv 0 \mod 4$.*

*Proof.* (1) If the parallel edges are positive, then they contain an extended Scharlemann cycle, which is impossible by Lemma 2.2(6). So assume that the edges are negative. It was shown in [Go, Lemma 4.2] that if the permutation is not transitive, then $M$ contains a cable space, which is impossible in our case because $M$ is hyperbolic. Hence the permutation is transitive. By the parity rule these edges connect parallel vertices in $\Gamma_2$, hence all the vertices of $\Gamma_2$ are parallel.

(2) and (3) are the same as [Wu1, Lemmas 1.7, 1.4, and 1.8]. □

**Lemma 2.4.** *Let $\widehat{e}$ be a family of parallel positive edges in $\Gamma_\alpha$. If the label $i$ appears twice among the labels of $\widehat{e}$, then there is a pair of edges among this family which form a Scharlemann cycle with $i$ as one of its labels. In particular, if $\widehat{e}$ has more than $n_\beta/2$ edges, then it contains a Scharlemann cycle.*

*Proof.* Since the edges are positive, by the parity rule $i$ cannot appear at both endpoints of a single edge in this family. Let $e_1, e_2, \ldots, e_k$ be consecutive edges of $\widehat{e}$ such that $e_1$ and $e_k$ have $i$ as a label. If $k \geq 4$, one can see that these edges contain an extended Scharlemann cycle in the middle, so by Lemma 2.2(6) we have $n_\beta = 2$, in which case $e_1, e_2$ is a Scharlemann cycle with $i$ as a label, as required. If $k = 2$ then $e_1, e_2$ is a Scharlemann cycle and we are done. If $k = 3$ the parity rule implies that the endpoints of $e_1$ and $e_3$ on the same vertex have the same label, so again $n_\beta = 2$ and $e_1, e_2$ is a Scharlemann cycle.

If $\widehat{e}$ has more than $n_\beta/2$ edges, then it has more than $n_\beta$ endpoints, so some label must appear twice. □

**Lemma 2.5.** (Properties of $\Gamma_2$) *(1) Any two Scharlemann cycles on $\Gamma_2$ have the same label pair.*

*(2) Any family of parallel positive edges $\widehat{e}$ in $\Gamma_2$ contains at most $n_1/2 + 1$ edges. Moreover, if it does contain $n_1/2 + 1$ edges, then two of them form a Scharlemann cycle.*

(3) Any family of parallel edges $\widehat{e}$ in $\Gamma_2$ contains at most $n_1$ edges.

(4) No three $i$-edges of $\Gamma_2$ are parallel for any $i$.

*Proof.* (1) This is [Wu3, Lemma 5.4(2)].

(2) If $\widehat{e}$ contains $n_1/2 + 1$ edges then they have $n_1 + 2$ endpoints, so some label $i$ appear twice. By Lemma 2.4 $i$ is the label of a Scharlemann cycle in this family. This proves the second sentence of (2).

Now suppose $\widehat{e}$ contains $n_1/2 + 2$ positive edges. If $n_1 = 2$ then there are only two nonparallel edges on the annulus $\widehat{F}_1$ connecting the two vertices, so a family of $n_1/2 + 2 = 3$ edges in $\widehat{e}$ would contain two which are parallel on both graphs, contradicting Lemma 2.2(2). If $n_1 > 2$, then 4 labels appear twice among the endpoints of the edges in $\widehat{e}$, so by Lemma 2.4 they are all labels of Scharlemann cycles in $\widehat{e}$, contradicting (1).

(3) If the edges in $\widehat{e}$ are positive, this follows from (2). So assume that $\widehat{e}$ consists of $n_1 + 1$ parallel negative edges $e_1, e_2, \ldots, e_{n_1+1}$. By Lemma 2.2(3) the first $n_1$ edges form several essential cycles in $\Gamma_1$. Let $e_1 \cup e_{j_2} \ldots \cup e_{j_k}$ be the one containing $e_1$. The edges $e_1$ and $e_{n_1+1}$ have the same label pair at their endpoints in $\Gamma_2$, hence have the same end vertices on $\Gamma_1$. By Lemma 2.2(2) they cannot be parallel on $\Gamma_1$, so after shrinking each vertex to a point, the path $e_{n_1+1}$ is homotopic rel $\partial$ to the path $e_{j_2} \cup \ldots \cup e_{j_k}$; in other words, the cycle $e_{j_2} \cup \ldots \cup e_{j_k} \cup e_{n_1+1}$ is an inessential cycle in $\Gamma_1$, which contradicts Lemma 2.2(3) because it is a cycle corresponding to the $n_1$ parallel edges $e_2, \ldots, e_{n_1+1}$ in $\Gamma_2$.

(4) If there are three $i$-edges which are parallel, then two of them would have the same label at one of their end vertices, so there would be at least $n_1 + 1$ parallel edges, contradicting (3). □

Let $u$ be a vertex of $\Gamma_\alpha$. Let $P, Q$ be two edge-endpoints on $\partial u$. Let $P_0 = P, P_1, \ldots, P_{k-1}, P_k = Q$ be the edge-endpoints encountered when traveling along $\partial u$ in the direction induced by the orientation of $u$. Then the *distance* from $P$ to $Q$ is defined as $\rho_u(P, Q) = k$. To emphasize that the distance is measured on $\Gamma_\alpha$, we may also use $\tau_\alpha(P, Q)$ to denote $\rho_u(P, Q)$. Notice that if the valency of $u$ is $m$, then $\rho_u(Q, P) = m - \rho_u(P, Q)$. The following lemma can be found in [Go].

**Lemma 2.6.** ([Go, Lemma 2.4]) *Let $a, b$ be components of $\partial F_\alpha \cap T_0$, and $x, y$ components of $\partial F_\beta \cap T_0$.*

*(i) Suppose $P, Q \in a \cap x$ and $R, S \in b \cap y$. If $\tau_\alpha(P, Q) = \tau_\alpha(R, S)$ then $\tau_\beta(P, Q) = \tau_\beta(R, S)$.*

*(ii) Suppose that $P \in a \cap x$, $Q \in a \cap y$, $R \in b \cap x$, and $S \in b \cap y$. If $\tau_\alpha(P, Q) = \tau_\alpha(R, S)$, then $\tau_\beta(P, R) = \tau_\beta(Q, S)$.* □

If $e_i$ is an edge with exactly one endpoint on $u$, $i = 1, 2$, define $\rho_u(e_1, e_2) = \rho_u(P_1, P_2)$, where $P_i = e_i \cap u$.

A pair of edges $e_1, e_2$ connecting two vertices $u, v$ in $\Gamma$ is an *equidistant pair* if $\rho_u(e_1, e_2) = \rho_v(e_2, e_1)$. If $u = v$, then define $e_1, e_2$ to be an equidistant pair if $\rho_u(P_1, P_2) = \rho_v(Q_2, Q_1)$, where $P_i, Q_i$ are the endpoints of $e_i$. One can check that this is independent of the choices of $P_i$ and $Q_i$.

The following lemma follows immediately from the definition, and is convenient for applications. In particular, it applies to parallel loops based at a vertex.

**Lemma 2.7.** *If $e_1, e_2$ is a pair of parallel edges connecting a pair of parallel vertices in $\Gamma$, then $e_1, e_2$ is an equidistant pair in $\Gamma$.* □

**Lemma 2.8.** *Let $e_1, e_2$ be a pair of edges with $\partial e_1 = \partial e_2$ in both $\Gamma_1$ and $\Gamma_2$. Then $e_1, e_2$ is an equidistant pair in $\Gamma_1$ if and only if it is an equidistant pair in $\Gamma_2$.*

*Proof.* Suppose $e_i$ joins vertices $a$ and $b$ in $\Gamma_1$, and vertices $x$ and $y$ in $\Gamma_2$, $i = 1, 2$. Denote the endpoints of $e_1, e_2$ in $\Gamma_1$ by $P, S$ and $Q, R$, respectively, as in Figure 2.1(a) or (c). There are two cases, depending on whether or not $P$ and $Q$ have the same label. We give a proof in the case that $P$ and $Q$ have the same label, say $x$. The proof for the other case is similar, using Lemma 2.6(2) instead of 2.6(1).

In this case we have $P, Q \in a \cap x$ and $R, S \in b \cap y$, so $e_1, e_2$ being equidistant on $\Gamma_1$ implies that
$$\tau_\alpha(P, Q) = \rho_a(e_1, e_2) = \rho_b(e_2, e_1) = \tau_\alpha(R, S).$$

By Lemma 2.6(1), we have $\tau_\beta(P, Q) = \tau_\beta(R, S)$, which implies that $\rho_x(e_1, e_2) = \rho_y(e_2, e_1)$, so $e_1, e_2$ is also equidistant on $\Gamma_2$. Similarly, $e_1, e_2$ being equidistant on $\Gamma_2$ also implies that they are equidistant on $\Gamma_1$. □

Figure 2.1

**Definition 2.9.** Given two slopes $r_1, r_2$ on $T_0$, choose a meridian-longitude pair $m, l$ on the torus $T_0$ so that $r_1 = m$, and the slope $r_2$ is represented by $dm + \Delta l$ for some $d$ with $1 \leq d \leq \Delta/2$. Then $d = d(r_1, r_2)$ is called the *jumping number* of $r_1, r_2$. Notice that if $\Delta = 4$, then $d = 1$, and if $\Delta = 5$, then $d = 1$ or 2.

**Lemma 2.10.** *Let $P_1, \ldots, P_\Delta$ be the points of $\partial u_i \cap \partial v_j$, labeled so that they appear successively on $u_i$. Let $d = d(r_1, r_2)$ be the jumping number of $r_1, r_2$. Then on $v_j$ these points appear in the order of $P_d, P_{2d}, \ldots, P_{\Delta d}$, in some direction. In particular, if $\Delta = 4$, then $P_1, P_2, P_3, P_4$ also appear successively on $v_j$.*

*Proof.* This follows by examining the order of the intersection points on the two curves $r_1$ and $r_2$. If $\Delta = 4$, we must have $d = 1$, hence the result follows. □

**Lemma 2.11.** ([HM, Proposition 5.1]) *If $C$ is a cycle in $\Gamma_\alpha$ consisting of positive $i$-edges, bounding a disk $D$ which contains no vertices in its interior, then $D$ contains a Scharlemann cycle.* □

**Lemma 2.12.** *(1) Suppose $G$ is a reduced graph on an annulus $A$ with $V$ vertices. Then $G$ has at most $3V - 2$ edges, and has a vertex of valency at most 5.*
*(2) A reduced graph $G$ on a torus has at most $3V$ edges.*

*Proof.* (1) By embedding the annulus $A$ in a sphere $S$, we can consider $G$ as a graph on $S$. Let $E, F$ be the number of edges and disk faces of $G$. Then $V - E + F \geq 2$ (the inequality may be strict if there are some non-disk faces.) All but two faces of $G$ have at least three edges, and each of the two exceptional ones has at least one edge. Hence we have
$$2 + 3(F - 2) \leq 2E.$$
Solving those two inequalities, we get $E \leq 3V - 2$.

If each vertex of $G$ has valency at least 6, then $6V \leq 2E$, which is impossible because $E \leq 3V - 2$.

(2) The calculation for $G$ on torus is similar but simpler. We have $V - E + F \geq 0$, and $3F \leq 2E$. The result follows by solving these inequalities. □

## §3. Case 1: One of the graphs has a single vertex

From now on we assume that $\widehat{F}_1$ is an essential annulus $A$ in $M(r_1)$, and $\widehat{F}_2$ is an essential torus $T$ in $M(r_2)$. Recall that $M(r_\alpha)$ is assumed to be irreducible. By [Go, Theorem 1.3] we may, and will always, assume that $\Delta \leq 5$. In this section we will consider the case where one of the graphs has only one vertex.

The link in Figure 7.1(a) is called the *Whitehead link*.

**Lemma 3.1.** *Suppose $n_1 = 1$ and $\Delta \geq 4$. Then $n_2 = 2$, $\Delta = 4$, and the graphs $\Gamma_i$ are as shown in Figure 3.1. The manifold $M$ is the Whitehead link exterior.*

*Proof.* In this case, there is only one edge in the reduced graph $\widehat{\Gamma}_1$, so all edges of $\Gamma_1$ are parallel positive edges. By the parity rule, all edges of $\Gamma_2$ are negative. In particular, $n_2 \geq 2$. Moreover, the $\Delta n_2 / 2$ parallel edges of $\Gamma_1$ contains an extended Scharlemann cycle, so by Lemma 2.2(6) we must have $n_2 = 2$. In this case there are exactly $\Delta$ edges in each graph. Since they are parallel on $\Gamma_1$, by Lemma 2.2(2) no two of them can be parallel on $\Gamma_2$. There are at most 4 nonparallel edges on the torus connecting two vertices. Hence $\Delta \leq 4$, and if $\Delta = 4$, then up to homeomorphism of the torus $T$ the graph is that shown in Figure 3.1. We will postpone until Lemma 7.2 the proof that if the intersection graphs are as shown in Figure 3.1 then the manifold $M$ is the exterior of the Whitehead link. This will complete the proof of Lemma 3.1. □

Figure 3.1

**Lemma 3.2.** *Suppose $n_2 = 1$ and $\Delta \geq 4$. Then $n_1 = 2$, $\Delta = 4$, and the graphs $\Gamma_i$ are as shown in Figure 3.3(c) and (d). The manifold $M$ is the Whitehead link exterior.*

*Proof.* The reduced graph $\widehat{\Gamma}_2$ on the torus is a subgraph of that shown in Figure 3.2(a), which has only three edges. The parity rule implies that $n_1 > 1$. The graph $\Gamma_2$ contains $\Delta n_1/2 \geq 2n_1$ edges, and so has a family of more than $n_1/2$ parallel positive edges. By Lemmas 2.5(2) and 2.2(4) $n_1$ must be even.

Figure 3.2

We want to show that $n_1 = 2$. If $\Delta = 5$ and $n_1 \geq 4$, then $\Gamma_2$ contains a family of $\Delta n_1/6 = \frac{1}{2}n_1 + \frac{1}{3}n_1 > \frac{1}{2}n_1 + 1$ parallel positive edges, which contradicts Lemma 2.5(2). Assume $\Delta = 4$. Then the labels at a family of parallel arcs are as shown in Figure 3.2(b). In particular, the number of parallel edges is even by the parity rule. Thus if $n_1 \geq 4$, then there are at least 4 edges in one of the families. From the labeling in Figure 3.2 we see that this family contains an extended Scharlemann cycle, contradicting Lemma 2.2(6).

Now we assume that $n_1 = 2$. In this case there are only two edges in $\widehat{\Gamma}_1$. Since no two edges can be parallel on both graphs, there are at most two edges in each family of parallel edges in $\Gamma_2$. Thus when $\Delta = 5$, there is only one possible configuration, shown in Figure 3.3(a), and when $\Delta = 4$, there are two possibilities, shown in Figure 3.3(b) – (c). The first two cases can be ruled out because there is an edge with the same label on both endpoints, which contradicts the parity rule. It follows that Figure 3.3(c) is the only possibility, which also uniquely determines the graph $\Gamma_1$, shown in Figure 3.3(d).

Figure 3.3

One can build the manifold $M$ as in Lemma 7.2 to show that $M$ is the exterior of the Whitehead link. However, there is a shortcut. Notice that the union of $A \times I$ with the Dehn filling solid torus $V_1$ is a genus 3 handlebody. After attaching the 2-handles coming from the faces of $\Gamma_2$ and then capping off the spherical boundary components, we get a manifold $Y$ with boundary a torus. Since the original manifold $M$ is assumed to be hyperbolic, this implies that $Y$ must be homeomorphic to $M(\gamma_1)$; in other words, there cannot be anything else attached to $Y$. Therefore $M$ is uniquely determined by the intersection graphs. Combining this with the above result, we conclude that there is at most one manifold $M$ such that $M(\gamma_1)$ contains a torus intersecting the Dehn filling solid torus only once, $M(\gamma_2)$ is annular, and $\Delta(\gamma_1, \gamma_2) > 3$. On the other hand, the Whitehead link exterior does satisfy those conditions, because the $-4$ surgery is annular, and the Seifert surface of one component of the link extends, after 0 surgery, to a torus which intersects the Dehn filling solid torus just once. Therefore, this must be the manifold. □

**Definition 3.3.** The graph $\Gamma_1$ is *k-separable* if there is an annulus $A'$ on $A$, such that (i) $A'$ contains exactly $k$ vertices of $\Gamma_1$, (ii) $\partial A'$ consists of two essential curves in the interior of $A$, and (iii) $\partial A'$ is disjoint from $\Gamma_1$. Denote by $\Gamma_1'$ the graph $\Gamma_1 \cap A'$ on $A'$, and by $\Gamma_2'$ the subgraph of $\Gamma_2$ consisting of edges which, when considered as edges in $\Gamma_1$, lie in $A'$. Notice that by definition $\Gamma_1$ is always $n_1$-separable.

**Lemma 3.4.** *Suppose $\Delta \geq 4$, and $\Gamma_1$ is 1-separable. Then $n_1 = 1$.*

*Proof.* Let $\Gamma_i'$ be the graphs on $A'$ and $T$ defined in Definition 3.3. As in the proof of Lemma 3.1, the graphs $\Gamma_i'$ must be the ones shown in Figure 3.1. Using the construction in the proof of Lemma 7.2, we see that the manifold $W = (A' \times I) \cup J_1 \cup (D \times I)$ is the Whitehead link exterior, where $D$ is a disk face of $\Gamma_2'$. Let $T_1 = \partial M - T_0$. Then $T_1$ contains the boundary of $A$, so it is nonempty. Let $T'$ be the component of $\partial W$ in the interior of $M$. Note that $T'$ separates $T_0$ from $T_1$, so $T'$ is incompressible in $M$ because otherwise $M$ would be reducible. Since $M$ is hyperbolic, $T'$ must be parallel to a boundary component of $M$. Since $W$ is not a product, we see that $T'$ must be parallel to $T_1$. In particular, $T_1$ is a torus. Now the annulus $A'$ can be extended to an essential annulus in $M(r_1)$ intersecting $J_1$ just once. By the minimality of $n_1$, we must have $n_1 = 1$. □

**Proposition 3.5.** *Suppose $\Delta \geq 4$. If $n_\alpha = 1$, then $\Delta = 4$, $n_\beta = 2$, and the manifold $M$ is the exterior of the Whitehead link shown in Figure 7.1.*

*Proof.* This follows immediately from Lemmas 3.1 and 3.2. □

§4. THE CASE THAT $n_2 = 2$

So far we have seen that if $\Delta = 4$ and one of the $n_\alpha = 1$ then the manifold $M$ is the exterior of the Whitehead link. The goal of this section is to show that if $n_2 = 2$ then the graph $\Gamma_1$ is 2-separable; see Proposition 4.6 below. In the next section it will be shown that if $\Gamma_1$ is 2-separable then $M$ must be the exterior of one of the links in Figure 7.1.

**Lemma 4.1.** *Suppose that $\Delta \geq 4$. If all vertices of $\Gamma_2$ are parallel, then $n_1 = n_2 = 2$, and hence $\Gamma_1$ is 2-separable.*

*Proof.* By Lemma 2.12 the reduced graph $\widehat{\Gamma}_2$ has at most $3n_1$ edges. For each $i$ there are $\Delta n_1$ $i$-edges in $\Gamma_2$, so two of them are parallel. By Lemma 2.4 this implies that $i$ is the index of some Scharlemann cycle in $\Gamma_2$. This is true for any index from 1 to $n_1$, but by Lemma 2.5(1) there are at most two labels that can appear as the labels of Scharlemann cycles in $\Gamma_2$, hence we have $n_1 = 2$. (We cannot have $n_1 = 1$ because of the parity rule.)

Now on $\Gamma_1$ there are only two families of parallel edges. Consider the function

$$\varphi : \{1, 2, \ldots n_2\} \to \{1, 2, \ldots n_2\}$$

defined so that if $e$ is an edge in $\Gamma_1$ with label $i$ on $u_1$ then it has label $\varphi(i)$ on the vertex $u_2$. Clearly this is well defined. Note that some family contains at least $2n_2$ parallel edges, so by Lemma 2.3(1) either $n_2 = 2$ and we are done, or the function $\varphi$ is transitive, i.e. it has only one orbit. In particular, if $n_2 > 2$ and $\varphi(k) = l$, then $\varphi(l) \neq k$.

On $\Gamma_2$ there are $\Delta n_1 > 3n_1$ edges, so there is a pair of adjacent parallel edges $e_1, e_2$ connecting, say, the vertex $v_k$ to $v_l$. Note that $e_1$ and $e_2$ have different labels at their endpoints on $v_k$, hence in $\Gamma_1$ the edge $e_1$ has label $k$, say, at the vertex $u_1$, and label $l$ at the vertex $u_2$, while $e_2$ has label $l$ at $u_1$ and label $k$ at $u_2$. Therefore $\varphi(k) = l$ and $\varphi(l) = k$, a contradiction. □

**Lemma 4.2.** *If $n_2 = 2$ then, up to homeomorphism of $T$, the reduced graph $\widehat{\Gamma}_2$ is a subgraph of the graph shown in Figure 4.1.*

*Proof.* If $\widehat{\Gamma}_2$ contains two loops based at $v_1$, then they cut the torus into a disk, hence there would be no loops based at $v_2$, so any edge with an endpoint on $v_2$ must have another endpoint on $v_1$; but since $v_1$ and $v_2$ have the same valency in $\Gamma_2$, this is impossible. Therefore $\widehat{\Gamma}_2$ has at most one loop based at each $v_i$. It is now easy to see that $\widehat{\Gamma}_2$ is isomorphic to a subgraph of that shown in Figure 4.1(a). □

Figure 4.1

We will use Figure 4.1(b) to indicate the reduced graph shown in Figure 4.1(a). The graph $\Gamma_2$ is obtained by replacing the edges in Figure 4.1(b) by families of parallel edges. Denote by $p_i$ the number of parallel edges in $\Gamma_2$ represented by $\widehat{e}_i$ in the reduced graph $\widehat{\Gamma}_2$.

**Lemma 4.3.** *If $\Delta \geq 4$ and $n_2 = 2$, then the two vertices $v_1, v_2$ of $\Gamma_2$ are antiparallel.*

*Proof.* Assume that $v_1, v_2$ are parallel. Then by Lemma 4.2 we have $n_1 = 2$. Since all edges of $\widehat{\Gamma}_2$ are positive, by the parity rule there is no loop on $\widehat{\Gamma}_1$, so $\widehat{\Gamma}_1$ consists of two families of parallel edges connecting $u_1$ to $u_2$.

Notice that in this case either all edges of $\widehat{\Gamma}_1$ have the same label on both endpoints, or they all have distinct labels on their two endpoints. In the first case $\widehat{\Gamma}_2$ consists of two families of parallel loops, so some pair of edges would be parallel on both graphs, which is impossible by Lemma 2.2(2). In the second case there are no loops in $\widehat{\Gamma}_2$, so there are only 4 edges in $\widehat{\Gamma}_2$. It follows again from Lemma 2.2(2) that $\Delta = 4$, and $\widehat{\Gamma}_2$ consists of 4 pairs of parallel edges. See Figure 4.2.

Consider a pair of parallel edges $e_1, e_2$ in $\Gamma_2$. Since $e_1, e_2$ are a pair of parallel edges connecting parallel vertices, by Lemma 2.7 they form an equidistant pair in $\Gamma_2$, so by Lemma 2.8 they also form an equidistant pair in $\Gamma_1$, that is $\rho_{u_1}(e_1, e_2) = \rho_{u_2}(e_2, e_1)$. Also, from Figure 4.2(b) we can see that $\rho_{u_2}(e_2, e_1) = 8 - \rho_{u_1}(e_1, e_2)$, so we must have $\rho_{u_1}(e_1, e_2) = \rho_{u_2}(e_2, e_1) = 4$. On the other hand, observe that the endpoints of $e_1, e_2$ on $u_1$ have different labels, so $\rho_{u_1}(e_1, e_2)$ must be an odd number, which is a contradiction. $\square$

Figure 4.2

**Lemma 4.4.** *Suppose $\Delta \geq 4$ and $n_2 = 2$. If $p_2 + p_3 \geq p_4 + p_5$, then $p_2 = p_3 = n_1$.*

*Proof.* By Lemma 2.5(3), $p_i \leq n_1$, $1 \leq i \leq 6$. Since $v_1$ and $v_2$ are antiparallel, we shall write $v_1 = v_+$, and $v_2 = v_-$. After relabeling if necessary, we may assume that the labels of the endpoints of the edges in $\widehat{e}_1$ are $1, 2, \ldots, k$ at the bottom of $v_+$. If $p_2 + p_3 = 2n_1 - r$, then since

$$\Delta n_1 = 2k + (p_2 + p_3) + (p_4 + p_5) \geq 2k + 2(2n - r),$$

we have $r \leq k$, so the labels of the left $r$ endpoints of edges in $\widehat{e}_1$ are $1, 2, \ldots, r$ at the top of $v_+$; see Figure 4.3(a). If $r > 2$, then the left $r$ edges of $\widehat{e}_1$ would contain an extended Scharlemann cycle, which would contradict Lemma 2.2(6) because $n_1 \geq k \geq r > 2$. Therefore we have $p_2 + p_3 \geq 2n_1 - 2$. Also, the left edge of $\widehat{e}_1$ cannot have the same label at both endpoints because of the parity rule. Hence either $p_2 + p_3 = 2n_1$, which implies $p_2 = p_3 = n_1$, or $p_2 + p_3 = 2n_1 - 2$. We need to rule out the second possibility.

Assume $p_2 + p_3 = 2n_1 - 2$. Since it is assumed that $p_2 + p_3 \geq p_4 + p_5$, we have

$$p_1 = p_6 = \frac{1}{2}(\Delta n_1 - (p_2 + p_3) - (p_4 + p_5)) \geq 2.$$

Examining the labels, we see that the two leftmost edges in $\widehat{e}_1$ form a Scharlemann cycle. Similarly for the two leftmost edges in $\widehat{e}_6$. By Lemma 2.5(1), these two Scharlemann cycles must have the same label pairs, say $\{1, 2\}$; see Figure 4.3(b). The existence of Scharlemann cycles implies that $n_1$ is even, and all vertices in $\Gamma_1$ of the same parity are parallel.

By Lemma 2.5(2) $\widehat{e}_1$ contains at most $n_1/2 + 1$ edges, so $p_4 + p_5 = \Delta n_1 - (p_2 + p_3) - 2p_1 \geq \Delta n_1 - 3n_1$. Since $p_4 + p_5 \leq p_2 + p_3$, we must have $\Delta = 4$.

First assume that $n_1 \geq 4$. In this case, because of the orientation of $v_+$ and $v_-$, the labels of the Scharlemann cycles must be as shown in Figure 4.3(b). It is easy to see that we cannot have $p_2 = p_3 = n_1 - 1$, otherwise an edge $e$ in the family $\widehat{e}_2$ would have labels of distinct parities on its two endpoints, which would contradict the parity rule because $e$ would connect antiparallel vertices in both graphs. Hence one of $p_2, p_3$ is $n_1 - 2$ and the other is $n_1$. Without loss of generality we may assume that $p_2 = n_1 - 2$, and $p_3 = n_1$. See Figure 4.3(b).

Let $e_1, \ldots, e_{n_1}$ be the parallel edges in the family $\widehat{e}_3$, where $e_i$ has label $i$ at its endpoint on $v_+$. From Figure 4.3(b) we can see that the other endpoint of $e_i$ is labeled $i+2$, hence these edges form two cycles $e_1 \cup e_3 \cup e_5 \cup \ldots \cup e_{n_1-1}$ and $e_2 \cup e_4 \cup e_6 \cup \ldots \cup e_{n_1}$ on $\Gamma_1$. Also, the edges of the Scharlemann cycle in the family $\widehat{e}_1$ form a cycle $e'_1 \cup e'_2$ on $\Gamma_1$ with vertices $u_1$ and $u_2$. By Lemmas 2.2(3) and 2.2(5) all these three cycles are essential on the annulus $A$, which implies that the cycle $e'_1 \cup e'_2$ must lie between the other two cycles. See Figure 4.3(c) for a possible picture of $\widehat{\Gamma}_1$, where $n_1 = 6$.

There are $\Delta = 4$ edge endpoints at $v_+$ labeled 1, among which the two endpoints on $e'_1$ and $e'_2$ are non-adjacent. Since $p_1 \leq n_1$, the other two edges labeled 1 at $v_+$ connect antiparallel vertices in $\Gamma_2$, hence connect parallel vertices in $\Gamma_1$. Since all the vertices with the same parity as $u_1$ are on one side of the cycle $e'_1 \cup e'_2$, this implies that all the edges with label $+$ at $u_1$ are on one side of $e'_1 \cup e'_2$. It follows that the endpoints of $e'_1$ and $e'_2$ on $u_1$ are adjacent among the endpoints of all such edges. Since we have seen that the corresponding endpoints of $e'_1$ and $e'_2$ are non-adjacent on $v_+$ among all endpoints labeled 1, this contradicts Lemma 2.10.

Now assume $n_1 = 2$. Since $p_1 \leq n_1$, and since $p_2 + p_3 = 2n_1 - 2 = 2 \geq p_4 + p_5$, we must have $p_2 + p_3 = p_4 + p_5 = p_1 = 2$. The graph $\Gamma_1$ is shown in Figure 4.3(d). Now the same argument as above applies: Among the 4 points of $\partial v_+ \cap \partial u_1$, the two on $e'_1$ and $e_2'$ are non-adjacent on $\partial v_+$ but adjacent on $\partial u_1$, which contradicts Lemma 2.10, completing the proof of the lemma. □

Figure 4.3

**Lemma 4.5.** *Suppose $\Delta \geq 4$, $n_2 = 2$, and $n_1 \geq 2$. Then $p_1 = p_6 > 0$.*

*Proof.* Assume that $p_1 = p_6 = 0$. The reduced graph $\widehat{\Gamma}_2$ has 4 edges $\widehat{e}_2, \widehat{e}_3, \widehat{e}_4, \widehat{e}_5$, each representing at most $n_1$ edges of $\Gamma_2$. The valency of $v_+, v_-$ in $\Gamma_2$ is $\Delta n_1$, and $\Delta \geq 4$, so we must have $\Delta = 4$, and each $\widehat{e}_i$ represents exactly $n_1$ edges.

Consider the permutation $\varphi : \{1, \ldots, n_1\} \to \{1, \ldots, n_1\}$ determined by the parallel edges in $\widehat{e}_i$. Since each family has exactly $n_1$ edges, one can see that $\varphi$ is independent of the choice of $\widehat{e}_i$. Let $L$ be the length of an orbit of $\varphi$. The edges in $\widehat{e}_i$ form cycles in $\Gamma_1$, each having length $L$. Thus the reduced graph $\widehat{\Gamma}_1$ consists of several cycles of length $L$, and each edge of $\widehat{\Gamma}_1$ represents exactly 4 edges of $\Gamma_1$.

CASE 1. $L \geq 3$.

There are 4 edges in $\Gamma_1$ from $u_1$ to $u_{\varphi(1)}$, two of which have label $+$ at $u_1$. On the other hand, there are 4 edges in $\Gamma_2$ having label 1 at $v_+$ and label $\varphi(1)$ at $v_-$, which is a contradiction. See Figure 4.4(a)–(b) for the case that $L = n_1 = 3$ and $\varphi(1) = 2$.

Figure 4.4

CASE 2. $L = 2$.

Let $k = \varphi(1) = n_1/2$. There is a pair of edges $a, b$ in the family $\widehat{e}_3$, such that $a$ has label 1 at $v_+$ and label $k$ at $v_-$, while $b$ has label $k$ at $v_+$ and label 1 at $v_-$. Let $P, S, Q, R$ be the endpoints of $a, b$, as shown in Figure 4.4(c). Notice that the vertices $u_1$ and $u_k$ are parallel, so from Figure 4.4(d) we can see that $\rho_1(a, b) = \rho_k(b, a)$. Hence $a, b$ are an equidistant pair in $\Gamma_1$. By Lemma 2.8 $a, b$ are also an equidistant pair in $\Gamma_2$. However, this is impossible because $\rho_+(a, b) = k$ and $\rho_-(b, a) = 7k$.

CASE 3. $L = 1$.

In this case all edges of $\Gamma_1$ are loops, hence the reduced graph $\widehat{\Gamma}_1$ consists of $n_1$ disjoint loops. Since $n_1 > 1$, $\Gamma_1$ is 1-separable, which contradicts Lemma 3.4. □

**Proposition 4.6.** *Suppose $\Delta \geq 4$, $n_2 = 2$, and $n_1 > 2$. Then $\Gamma_1$ is 2-separable.*

*Proof.* By Lemmas 4.3–4.5 we may assume that $v_1 = v_+$ and $v_2 = v_-$ are antiparallel, $p_1 = p_6 > 0$, and $p_2 = p_3 = n_1$. We may assume without loss of generality that the labels of the edge endpoints are as shown in Figure 4.5(a).

The families of parallel edges in $\widehat{e}_2$ and $\widehat{e}_3$ determine the same permutation $\varphi$ on the set of labels $\{1, \ldots, n_1\}$. Let $L$ be the length of an orbit of $\varphi$. Let $k = \varphi(1)$. The argument splits into two cases: $L \geq 2$, and $L = 1$. We will show that the first case leads to a contradiction, while the second implies that $\Gamma_1$ is 2-separable.

CASE 1. $L \geq 2$.

The parallel edges in $\widehat{e}_2$ form cycles in $\Gamma_1$. There is a loop (namely $C$ in Figure 4.5(a)) in $\Gamma_2$ at $v_+$ with endpoints labeled 1 and $n_1$, so by the parity rule $u_1$ and $u_{n_1}$ are of opposite signs, hence 1 and $n_1$ belong to distinct orbits of $\varphi$. Let $\theta$ be the cycle containing $u_1$ and $u_k$, and let $\theta'$ be the one containing $u_{n_1}$ and $u_{k-1}$. By Lemma 2.2(3) all these cycles are essential in $A$. Since there is an edge $C$ in $\Gamma_2$ with endpoints labeled 1 and $n_1$ (see Figure 4.5(a)), which connects $u_1$ to $u_{n_1}$ in $\Gamma_2$, there are no vertices inside the annulus between $\theta$ and $\theta'$.

First assume that $\Delta = 4$. Let $A, B, C, D$ be the edges with an endpoint on $v_+$ labeled 1, and let $A, B, F, G$ be those with an endpoint on $v_-$ labeled $k$. See Figure 4.5(a). On $\Gamma_2$ these edges form a subgraph as shown in Figure 4.5(b). The edges $C$ and $G$ cut the annulus between $\theta$ and $\theta'$ into two disks $U_1$ and $U_2$, as shown in the figure.

The endpoints of the edges $A, B, C, D$ appear successively on $\partial v_+$ in $\Gamma_2$, so by Lemma 2.10 they also appear successively on $\partial u_1$ in $\Gamma_1$. Similarly, the endpoints of $A, B, F, G$ appear successively on $\partial u_k$. From Figure 4.5(b) one can see that exactly one of the edges $D$ and $F$, say $D$, lies inside of $D$; in particular, we must have $D \neq F$. On the other hand, in $U_1$ the only vertices which are parallel to $u_1$ are $u_k$ and $u_1$ itself, so $D$ must have its other endpoint on $u_k$. But this is a contradiction, because in $\Gamma_2$ the endpoint of $D$ on $v_-$ would then have label $k$, implying that $D = F$.

Now assume that $\Delta = 5$. Let $A, B, C, D, E$ be the edges with an endpoint on $v_+$ labeled 1, and let $A, B, F, H, G$ be those with an endpoint on $v_-$ labeled $k$, as shown in Figure 4.5(a). By Lemma 2.10 there are two possible orders in which the endpoints of these edges appear on $\partial u_1$ and $\partial u_k$ respectively. (i) The jumping number $d = d(r_1, r_2) = 1$. In this case the endpoints appear as $A, B, C, D, E$ on $\partial u_1$, and the proof is similar to that given above, noticing that $D, E$ are in $U_1$ if and only if $F$ and $H$ are outside of $U_1$. (ii) The jumping number $d = 2$, so the endpoints appear in the order $A, D, B, E, C$ on $\partial u_1$, and $A, H, B, G, F$ on $\partial u_k$. See Figure 4.5(c). One can see that the edges $E$ and $F$ lie in different disks $U_1, U_2$, while $D$ and $H$ lie between $A$ and $B$; hence the argument above can be applied to the edges $E$ and $F$ give a similar contradiction.

Figure 4.5

CASE 2. $L = 1$.

In this case all the edges in $\widehat{e}_2$ and $\widehat{e}_3$ are loops in $\Gamma_1$. Orient all the loops from label $+$ to label $-$. Let $e_i$ and $e'_i$ be the edges in $\widehat{e}_2, \widehat{e}_3$ which form loops based at $u_i$. There are two possibilities, depending on whether they have the same orientation; see

Figure 4.6. We claim that *in either case all the non-loop edges with an endpoint at $u_i$ are on one side of the loops.*

First assume $\Delta = 4$. If the two edges have different orientations, then since their + endpoints are adjacent among all endpoints on $\partial v_+$ which are labeled $i$, by Lemma 2.10 these two endpoints are also adjacent on $\partial u_i$ among all endpoints labeled $+$. Since the $+$ and $-$ labels alternate around $\partial u_i$, it follows that all the endpoints besides those of $e_i$ and $e'_i$ are on one side of the loops. If the two edges have the same orientation, then there is another edge between them. Since the $+$ and $-$ labels alternate, from Figure 4.6(b) we see that the number of labels on each side of the loops must be even. Therefore the remaining two endpoints must be on the same side of the loops.

Now assume $\Delta = 5$. We have $p_4 + p_5 + 2p_1 = 3n_1$. Since $p_4, p_5 \leq n_1$, and $p_1 \leq n_1/2 + 1$, we have only two possibilities: Either $p_4 = p_5 = n_1$ and $p_1 = n_1/2$, or $p_4 + p_5 = 2n_1 - 2$ and $p_1 = n_1/2 + 1$. In the first case there are four loops based at each vertex, so there are only two adjacent edge endpoints left, which must be on the same side of the loops. In the second case, the argument beginning at the second paragraph of the proof of Lemma 4.4, with $p_2, p_3$ replaced by $p_4, p_5$, applies (although here we don't have $p_4 + p_5 \geq p_2 + p_3$), to show that the reduced graph $\Gamma_1$ contains a subgraph as shown in Figure 4.3(c), (which illustrates the case $n_1 = 6$.) However, since $p_1 = n_1/2 + 1 \geq 3$, there is an edge in $\widehat{e}_1$ adjacent to the Scharlemann cycle, with endpoints labeled 3 and $n_1$. On $\Gamma_1$ this edge connects $u_3$ to $u_{n_1}$, which is impossible because $u_3$ and $u_{n_1}$ are on different sides of the loop $e'_1 \cup e'_2$ in Figure 4.3(c). This completes the proof of the claim.

Consider the reduced graph $\widehat{\Gamma}_1$. We have shown that each vertex is the base of a loop. Let $u', u''$ be two vertices that are connected by an edge. The above claim says that there are no edges connecting $u', u''$ to any other vertices. Thus there is an annulus $A'$ containing only the vertices $u', u''$, and with boundary disjoint from $\Gamma_1$. By definition $\Gamma_1$ is 2-separable. □

Figure 4.6

## §5. The case that $\Gamma_1$ is 2-separable

Assume that $n_1, n_2 \geq 2$, and $\Delta \geq 4$. In Section 4 we have shown that if $n_2 = 2$ then $\Gamma_1$ is 2-separable. Lemma 5.1 below shows that the converse is also true. The goal of this section is to show that if $\Gamma_1$ is 2-separable, (in particular, if one of the $n_\alpha = 2$), then $n_1 = n_2 = 2$, and the manifold $M$ is the exterior of one of the links in Figure 7.1(b) or (c). See Proposition 5.5.

Throughout this section we will assume that $\Gamma_1$ is 2-separable. We use $A', \Gamma_1', \Gamma_2'$ to denote the annulus and graphs given in Definition 3.3. Let $u_1, u_2$ be the vertices of $\Gamma_1'$. The labeling of $u_i$ in $\Gamma_1'$ may not be the same as that in $\Gamma_1$, however this does not affect the proof below. Notice also that a pair of edges is an equidistant pair in $\Gamma_i'$ if and only if it is an equidistant pair in $\Gamma_i$, hence Lemma 2.8 still applies to $\Gamma_i'$. Also, since $\Gamma_i'$ are subgraphs of $\Gamma_i$, most of the properties of $\Gamma_i$ still hold for $\Gamma_i'$. For example, there is no extended Scharlemann cycle in $\Gamma_1'$.

**Lemma 5.1.** *Suppose $\Delta \geq 4$. If $\Gamma_1$ is 2-separable, then $n_2 = 2$.*

*Proof.* Let $A'$, $\Gamma_i'$ be as in Definition 3.3. The reduced graph $\widehat{\Gamma}_1'$ is a subgraph of the graph shown in Figure 5.1. If some edge represents more than $n_2$ parallel edges of $\Gamma_1'$, then by Lemma 2.3(1) either $n_2 = 2$ and we are done, or all vertices of $\Gamma_2$ are parallel, which by Lemma 4.1 again implies that $n_2 = 2$. We can thus assume that $\Delta = 4$, and that each edge in Figure 5.1 represents exactly $n_2$ parallel edges of $\Gamma_1'$. The $n_2$ parallel loops based at $u_1$ contains a Scharlemann cycle, so by Lemma 2.2(4) $n_2$ is even. Also, notice that the labels of the $2n_2$ endpoints of these loops appear successively, so if $n_2 \geq 4$ then the 4 edges in the middle would form an extended Scharlemann cycle, which is impossible by Lemma 2.2(6). □

Figure 5.1

**Lemma 5.2.** *Suppose $\Delta \geq 4$, and $\Gamma_1$ is 2-separable. Then $\Gamma_2'$ must be one of the graphs shown in Figure 5.2(d) or (e).*

*Proof.* By Lemma 5.1 we have $n_2 = 2$, by Lemma 4.3 the two vertices of $\Gamma_2'$, $v_+$ and $v_-$, are antiparallel, and by Lemma 4.2 the reduced graph $\widehat{\Gamma}_2'$ is a subgraph of the graph shown in Figure 4.1. Denote by $p_i'$ the number of edges of $\Gamma_2'$ in the family $\widehat{e}_i$ shown in Figure 4.1(b). Recall that $p_i$ is the number of edges of $\Gamma_2$ in $\widehat{e}_i$. Without loss of generality we may assume that $p_2 + p_3 \geq p_4 + p_5$. Then by Lemma 2.5(3) we have $p_i \leq 2$, and by Lemma 4.4 we have $p_2 = p_3 = n_1$. By considering the edges of $\Gamma_2'$ in each family $\widehat{e}_i$, it follows that $p_i' \leq 2$, and $p_2' = p_3' = 2$. According as $p_1' = 0, 1$ or 2, and $\Delta = 4$ or 5, we have the following possibilities.

(a) $p_1' = 0$, $\Delta = 4$, $p_4' + p_5' = 4$;
(b) $p_1' = 1$, $\Delta = 4$, $p_4' + p_5' = 2$;
(c) $p_1' = 1$, $\Delta = 5$, $p_4' + p_5' = 4$;
(d) $p_1' = 2$, $\Delta = 4$, $p_4' + p_5' = 0$;
(e) $p_1' = 2$, $\Delta = 5$, $p_4' + p_5' = 2$;

Using the parity rule, each of the above uniquely determines $\Gamma_2'$ up to a homeomorphism of the underlying torus $T$. See Figure 5.2(a)–(e). We want to show that cases (a), (b) and (c) cannot occur.

Consider the graph in Figure 5.2(a). It is easy to see that either all edges have the same label at both endpoints, or each edge has distinct labels at its two endpoints. In the first case $\Gamma_1'$ consists of loops, so $\Gamma_1$ is 1-separable. By Lemma 3.4 we must then have $n_1 = 1$, which is a contradiction. In the second case, the two vertices of $\Gamma_1'$ are parallel by the parity rule, the reduced graph $\widehat{\Gamma}_1'$ consists of two edges from $u_1$ to $u_2$, so one can see that any pair of edges is an equidistant pair on $\Gamma_1'$. On the other hand a pair of parallel edges $\{e_1, e_2\}$ is not an equidistant pair on $\Gamma_2'$ because $\rho_{v_+}(e_1, e_2) = 1$ and $\rho_{v_-}(e_2, e_1) = 7$. This contradicts Lemma 2.8, showing that $\Gamma_2'$ cannot be the graph in Figure 5.2(a).

In all the remaining cases, $p_1' \geq 1$, so there is a loop based at $v_+$. By the parity rule, the two vertices of $\Gamma_1'$ are antiparallel. Thus any edge in $\Gamma_2'$ from $v_+$ to $v_-$ must have the same label at its two endpoints, and is a loop in $\Gamma_1'$. The reduced graph of $\Gamma_1'$ is a subgraph of that in Figure 5.1, so any two loops based at $u_1$ are parallel. In particular, the edges $e_1, e_2$ shown in Figure 5.2(a)–(b) are parallel loops in $\Gamma_1'$, so by Lemma 2.7 they form an equidistant pair in $\Gamma_1'$, and by Lemma 2.8 they must also be an equidistant pair in $\Gamma_2'$. However, in both cases we have $\rho_{v_+}(e_1, e_2) = 2$ and $\rho_{v_-}(e_2, e_1) = 4$. Hence cases (b) and (c) cannot happen. □

Figure 5.2

**Lemma 5.3.** *Each of the graphs $\Gamma'_2$ in Figure 5.2(d) and (e) uniquely determines the graph $\Gamma'_1$ up to symmetry, as shown in Figure 5.3(a) and (b), respectively. The correspondence between the edges in $\Gamma'_1$ and $\Gamma'_2$ is also uniquely determined up to symmetry.*

*Proof.* The reduced graph $\widehat{\Gamma}'_1$ is a subgraph of the graph shown in Figure 5.1. The edges of $\Gamma'_2$ come in parallel pairs, which are nonparallel in $\Gamma'_1$ by Lemma 2.2(2). Using the fact that any loop in $\Gamma'_2$ is a non-loop in $\Gamma'_1$, and vice versa, we see that the graphs $\Gamma'_1$ corresponding to $\Gamma'_2$ in Figure 5.2(d) and (e) are the ones in Figure 5.2(a) and (b), respectively. It remains to show that the correspondence between the edges of $\Gamma'_1$ and $\Gamma'_2$ is unique up to symmetry.

Label the edges of $\Gamma'_2$ as in Figure 5.2(d). Up to symmetry we may assume that the edge $f$ is the one shown in Figure 5.3(a). The points of $\partial v_+ \cap \partial u_1$ appear in the order of $e, f, a, c$ on $\partial v_+$, so by Lemma 2.10 they must also appears in this order on $\partial u_1$. This determines these edges in $\Gamma'_1$. Similarly, the points of $\partial v_+ \cap \partial u_2$ appear in the order of $f, e, b, d$ on $\partial v_+$, so they also appear in this order on $\partial u_2$. This determines the edges $b, d$. By looking at the points of $\partial v_- \cap \partial u_1$, one can further determine the edges $h, g$. Hence in case (d) the edge correspondence is determined by $\Gamma'_2$.

The proof of case (e) is similar, except that in this case there is a jumping number $d = d(r_1, r_2)$ given by Definition 2.9, so if $P_1, P_2, \ldots, P_5$ are the points of $\partial v_i \cap \partial u_j$, numbered succesively on $\partial v_i$, (here $i = +$ or $-$), then they appear as $P_d, P_{2d}, \ldots, P_{5d}$ on $\partial u_j$. Label the edges of $\Gamma'_2$ as in Figure 5.2(e), and fix a choice of $e$ on $\Gamma'_1$. Then there is only one more edge from $u_1$ to $u_2$ with endpoint labeled $+$ on $u_1$, so that must be the edge $f$. The points of $\partial v_+ \cap \partial u_1$ appear as $e, j, f, a, c$ on $\partial v_+$, and as $e, *, *, *, f$ on $\partial u_1$, so we must have $d = 2$, and the edges are in the order of $e, a, j, c, f$. Similarly one can determine the other edges of $\Gamma'_1$. See Figure 5.3(b). □

Figure 5.3

Let $F'_1$ be the punctured annulus $A' - \text{Int}(u_1 \cup u_2)$. Let $F'_2$ be the subsurface of $F_2$ consisting of disk faces of $\Gamma'_2$. Consider the submanifold $X$ of $M$ which is a regular neighborhood of $T_0 \cup F'_1 \cup F'_2$ with spherical boundary components capped off by 3-balls. Since $M$ is irreducible, each sphere on $\partial N(T_0 \cup F'_1 \cup F'_2)$ does bound a 3-ball in $M$, hence $X$ is indeed a submanifold of $M$. The manifold $X$ can be constructed as follows. Lut $Y = (F'_1 \times I) \cup (T_0 \times I)$. For each disk face $D_i$ of $\Gamma'_2$ attach a 2-handle $H_i = D_i \times I$ to $Y$ according to the boundary curve of $D_i$ on $F'_1 \cup T_0$, and finally cap off the spherical boundaries.

**Lemma 5.4.** *The manifold $X$ constructed above is uniquely determined by the graph $\Gamma'_2$.*

*Proof.* By Lemma 5.3 the graph $\Gamma'_1$ and the correspondence between the edges of $\Gamma'_1$ and $\Gamma'_2$ are uniquely determined by $\Gamma'_2$. Thus for each face $D_i$ of $\Gamma'_2$, the boundary curve of $D_i$ on $F'_1 \cup T_0$ is determined. The lemma now follows from the above construction of $X$, because $X$ is uniquely determined by the way the 2-handles $H_i$ are attached to $Y$. □

**Proposition 5.5.** *Suppose $\Delta \geq 4$, and $\Gamma_1$ is 2-separable. Then $n_1 = n_2 = 2$, and there are exactly two such manifolds $M$. If $\Delta = 4$, $M$ is the exterior of the link in Figure 7.1(b), and if $\Delta = 5$, $M$ is the exterior of the link in Figure 7.1(c).*

*Proof.* By Lemma 5.1 we have $n_2 = 2$. Let $F_1', X, Y$ be as above. The surface $F_1'$ separates $Y$ into two sides, the black side $Y_1$ and the white side $Y_2$. Let $P_i$ be the punctured torus $(\partial Y - T_0) \cap Y_i$. Each face of $\Gamma_2'$ is colored according to the side it is on. Since $F_2$ is transverse to $F_1'$, the colors on the faces of $\Gamma_2'$ alternate. From Figure 6.2(d) and (e) we see that all the white disk faces have three edges, so they are nonseparating curves on the punctured torus $P_2$. Thus after attaching the 2-handles and capping off the spheres, the surface $P_2$ becomes an annulus. On the black side, there is a Scharlemann cycle on $\Gamma_2'$, so after attaching the $H_i$ and capping off the spheres, $P_1$ become either an annulus or a pair of disks. The second case does not happen because otherwise the boundary of $X$ contains a sphere separating the boundary components of $M$, which contradicts the fact that $M$ is irreducible. It follows that $X$ has exactly two boundary components, each being a torus. The annulus $A'$ is an essential annulus in $X$, with one boundary component on the outside torus $T'$, and two other boundary components on $T_0$, hence $X$ is not a product $T_0 \times I$.

The torus $T'$ is incompressible in the global manifold $M$, because otherwise after compression it would become a sphere separating the boundary components of $M$, contradicting the irreducibility of $M$. Since $M$ is atoroidal, and $X$ is not a product, it follows that $T'$ must be parallel to another boundary components $T_1$ of $M$. That is, $M - \text{Int} X$ is a product $T_1 \times I$. In particular, by Lemma 5.4 $M \cong X$ is determined by $\Gamma_2'$. Since there is only one possible $\Gamma_2'$ for each of $\Delta = 4, 5$, given in Figure 5.2(d) and (e), we see that there is at most one such manifold $M$ for each of $\Delta = 4, 5$. The existence of such manifolds will be established by Theorem 7.5.

Since $M - \text{Int} X \cong T_1 \times I$, the annulus $A'$ in $X(r_1)$ can be extended to an essential annulus $A''$ in $M(r_1)$ intersecting $J_1$ just twice. By the minimality of $n_1$, we must have $n_1 = 2$. □

§6. THE GENERIC CASE

After the results in the previous sections, we may now assume that $n_\alpha \geq 3$ for $\alpha = 1, 2$. The purpose of this section is to show that in this case we must have $\Delta \leq 3$.

**Lemma 6.1.** *If $n_1 \geq 3$, and $n_2 = 3$, then $\Delta \leq 3$.*

*Proof.* By Lemma 4.2 this is true if all vertices of $\Gamma_2$ are parallel, so we may assume that $v_1, v_2$ are positive vertices, and $v_3$ is a negative vertex of $\Gamma_2$. Assume $\Delta \geq 4$.

First assume that there is only one family of parallel loops based at $v_3$, containing $k$ edges. By Lemma 2.5(3), $k \leq n_1$. These loops contribute $2k$ endpoints at $v_3$, so there are $\Delta n_1 - 2k$ negative edges at $v_3$, hence the number of 3-edges in $\Gamma_1$ is

$$E = (\Delta n_1 - 2k) + k = \Delta n_1 - k \geq 3n_1.$$

By Lemma 2.2(2) no two negative 3-edges are parallel in $\Gamma_1$ because they are parallel positive edges in $\Gamma_2$. Since $n_2 = 3$, $T$ is nonseparating, so by Lemmas 2.4 and 2.2(4) no two positive 3-edges are parallel in $\Gamma_1$. Thus all the 3-edges are mutually nonparallel, so the reduced graph $\widehat{\Gamma}_1$ has at least $3n_1$ edges, which contradicts Lemma 2.12(1).

Now assume that there are at least two families of parallel loops based at $v_3$. These (together with $v_3$) cut $T$ into disks, so there are no loops based at the other vertices, and there is at most one family of parallel edges from $v_1$ to $v_2$. The argument in this case is now similar to that given above, considering the 1-edges of $\Gamma_1$ instead of the 3-edges. □

Let $G$ be a graph on a disk $D$. A vertex $v$ of $G$ is an *interior vertex* if there is no arc from $v$ to $\partial D$ with interior disjoint from $G$. It is a *cut vertex* if there is an arc $\alpha$ on $D$, cutting $D$ into $D_1, D_2$, such that the graphs $G_i = D_i \cap G$ have the properties that $G_1 \cap G_2 = v$, and neither $G_i$ is a single vertex. In this case write $G = G_1 *_v G_2$. A vertex $v$ is a *boundary vertex* if it is neither an interior vertex nor a cut vertex. Define

$a_i = a_i(G)$: the number of boundary vertices $v$ with $val(v, G) = i$;
$l = l(G)$: the number of trivial loops in $G$;
$p = p(G)$: the number of pairs of adjacent parallel edges;
$\sigma(G) = 3a_1 + 2a_2 + a_3 + 4l + 2p$;
$\tau(G) = a_1 + a_2 + l + p$.

**Lemma 6.2.** *Let $G$ be a connected graph on a disk $D$, which is not a single vertex. Suppose each interior vertex $v$ of $G$ has $val(v) \geq 6$. Then (i) $\sigma(G) \geq 6$, and (ii) if $G$ has no interior vertices then $\tau(G) \geq 2$.*

*Proof.* The result is obvious if $G$ has only one edge. Also, if $G = G_1 *_v G_2$ then the vertex $v$ contributes at most 3 to $\sigma(G_i)$ and at most 1 to $\tau(G_i)$, so

$$\sigma(G) \geq (\sigma(G_1) - 3) + (\sigma(G_2) - 3) \geq \sigma(G_1) + \sigma(G_2) - 6;$$
$$\tau(G) \geq (\tau(G_1) - 1) + (\tau(G_2) - 1) \geq \tau(G_1) + \tau(G_2) - 2;$$

hence the result follows by induction. We can thus assume that $G$ contains more than one edge, and has no cut vertices. By considering a sub-disk if necessary, we may assume that $\partial D \subset G$. Notice that in this case $a_1(G) = 0$.

(i) Taking the double of $D$, we get a graph $G'$ on a sphere $S$, with $a_0(G') = a_1(G') = a_3(G') = a_5(G') = 0$, $a_2(G') = a_2(G)$, and $a_4(G') = a_3(G)$. There are $2l$ faces with 1 edge, $2p$ faces with 2 edges, and all the others have at least 3 edges. Denote by $V, E, F$ the number of vertices, edges, and faces of $G'$. We have the following inequalities:

$$2a_2 + 4a_3 + 6(V - a_2 - a_3) \leq 2E$$
$$2l + 2(2p) + 3(F - 2l - 2p) \leq 2E$$
$$V - E + F = 2$$

Solving these inequalities gives $\sigma(G) = 2a_2 + a_3 + 4l + 2p \geq 6$.

(ii) This is obvious if all edges of $G$ are on the boundary of $D$. So assume that $e$ is an edge which is not on the boundary of $D$. It cuts $D$ into two disks $D_1, D_2$, and $G$ into $G_1, G_2$. If $l(G_i) + p(G_i) \geq 1$ for both $i$, then $\tau(G) \geq 2$ and we are done. Suppose $l(G_i) = p(G_i) = 0$ for some $i$. Then by induction we have $\tau(G_i) = a_2(G_i) \geq 2$, i.e. it has at least two boundary vertices of valency 2. If both endpoints of $e$ have valency 2 in $G_i$, then shrinking $e$ to a point gives a graph $G'_i$, which by induction must have another boundary vertex of valency 2. Hence $G_i$ must have a boundary vertex of valency 2 which is not on $\partial e$. Similarly for the other graph $G_j$. Hence we have $\tau(G) \geq 2$. □

**Lemma 6.3.** *Suppose $n_1 \geq 3$, $n_2 = 4$, and $\Delta \geq 4$. Then $\Gamma_1$ cannot have 4 length 2 Scharlemann cycles on mutually distinct label pairs.*

*Proof.* There are only four possible label pairs $\{i, i+1\}$, $i = 1, \ldots, 4$, so if the lemma were false then they would all appear as label pairs of Scharlemann cycles in $\Gamma_1$. Denote by $V_i$ the part of the Dehn filling solid torus between $v_i$ and $v_{i+1}$. The existence of Scharlemann cycles in $\Gamma_1$ implies that the torus $T$ in $M(r_2)$ is separating (Lemma 2.2(4)). Without loss of generality we may assume that $V_1, V_3$ lie on the side containing some boundary components of $M(r_2)$.

Let $\mathcal{E}^i = \{e_1^i, e_2^i\}$ be the Scharlemann cycle on $\Gamma_1$ with label pair $\{i, i+1\}$, and let $D_i$ be the disk bounded by $\mathcal{E}^i$ on the punctured annulus $F^1$. Then $e_1^1 \cup e_2^1$ is a cycle on $\Gamma_2$ with endpoints on $v_1, v_2$. Consider the manifold $X = (T \times I) \cup V_1 \cup (D_1 \times I)$. It is easy to see that $X$ is a cable space, whose boundary consists of two tori $T \cup T'$, with $T' \cap J_2$ consisting of two disks. If $T'$ is compressible, then after compression it becomes a sphere containing $T$ on one side and some boundary component of $M$ on the other, which is therefore a reducing sphere in $M(r_2)$; by [Wu3] in this case we would have $\Delta \leq 2$. So assume that $T'$ is incompressible. Then $X' = (T' \times I) \cup V_3 \cup (D_3 \times I)$ is also a cable space, hence $T'$ is not parallel to a torus in $\partial M(r_2)$. It follows that $T'$ is an essential torus in $M(r_2)$ intersecting $J_2$ in fewer than $n_2$ disks, contradicting the choice of the essential torus $T$. □

A vertex $u$ in $\Gamma_1$ is a *full vertex* if its valency in $\Gamma_1^+$ is $\Delta n_2$. In other words, $u$ is full if all the edges incident to it are positive.

**Lemma 6.4.** *Suppose $n_1 \geq 3$, $n_2 \geq 4$, and $\Delta \geq 4$. Then each full vertex $u_j$ of $\Gamma_1$ has valency at least 6 in $\widehat{\Gamma}_1$.*

*Proof.* Assume that $val(u_j, \widehat{\Gamma}_1) \leq 5$. There are a total of $\Delta n_2$ edge endpoints at $u_j$, and each family of parallel edges contains at most $n_2/2 + 2$ edges (Lemma 2.3(3)), so we have
$$\frac{\Delta n_2}{5} \leq \frac{n_2}{2} + 2.$$

This immediately gives $n_2 \leq 6$. If $n_2 = 6$, there is a family of 5 parallel edges, with endpoints labeled $1, 2, 3, 4, 5$ at $u_j$, say. Then one can see that this family contains two Scharlemann cycles on label pairs $\{1, 2\}$ and $\{4, 5\}$, respectively, which contradicts Lemma 2.3(2). If $n_2 = 5$ then $T$ is non-separating, so by Lemmas 2.4 and 2.2(4) each family of parallel edges contains at most $n_2/2$ edges. Replacing the right hand side of the above inequality with $n_2/2$, we see that $\Delta < 3$, contradicting the assumption. Hence we must have $n_2 = 4$.

Consider the reduced negative graph $\widehat{\Gamma}_2^-$. Let $E, F$ be the number of its edges and disk faces. There are four edges in $\Gamma_2^-$ labeled $j$ at each $v_i$, no two of which are parallel because otherwise there would be $n_1 + 1$ parallel edges in $\Gamma_2$, which would contradict Lemma 2.5(3). Thus the valency of $v_i$ in $\widehat{\Gamma}_2^-$ is at least 4 for any $i$, so

$$4n_2 \leq 2E.$$

Each face of $\widehat{\Gamma}_2^-$ must have at least 4 edges because two adjacent vertices on the boundary of the face have opposite signs; hence

$$4F \leq 2E.$$

Calculating the Euler number we have

$$n_2 - E + F \geq 0.$$

These three inequalities hold if and only if they are all equalities. In particular, each vertex has valency 4 in $\widehat{\Gamma}_2^-$, and there are exactly $2n_2 = 8$ edges in $\widehat{\Gamma}_2^-$.

The four edges of $\widehat{\Gamma}_2^-$ incident to $v_2$ have the other endpoint on $v_1$ or $v_3$. If there were three edges of $\widehat{\Gamma}_2^-$ connecting $v_2$ to $v_1$, then they would cut the torus $T$ into a disk. Since $val(v_3, \widehat{\Gamma}_2^-) = 4$, some pair of edges of $\widehat{\Gamma}_2^-$ incident to $v_3$ would be parallel, which is absurd because $\widehat{\Gamma}_2^-$ is a reduced graph. Therefore there are exactly two edges of $\widehat{\Gamma}_2^-$ with endpoints on $v_1$ and $v_2$. Similarly, for each $i = 1, \ldots, 4$, there are exactly two edges in $\widehat{\Gamma}_2^-$ connecting $v_i$ to $v_{i+1}$.

The 8 edges of $\widehat{\Gamma}_2^-$ contain the 16 $j$-edges, and by Lemma 2.5(4) each edge of $\widehat{\Gamma}_2^-$ contains at most two $j$-edges, and hence exactly two $j$-edges. Combining this with the conclusion of the last paragraph, we see that there are 4 $j$-edges connecting $v_i$ to $v_{i+1}$ for any $i$. On $\Gamma_1$ this means that among the 16 edges incident to $u_j$, exactly 4 of them have label pairs $\{i, i+1\}$, for each $i$.

By Lemma 6.3, some pair $\{i, i+1\}$, say $\{1, 2\}$, is not the label pair of a length 2 Scharlemann cycle on $\Gamma_1$. Thus by Lemma 2.4 the 4 edges with label pair $\{1, 2\}$ must be mutually nonparallel. Notice that if an edge $e$ is parallel to an edge with label pair $\{1, 2\}$, then it has label pair either $\{1, 2\}$ or $\{3, 4\}$. Since there are only 5 families of parallel edges at $u_j$, 4 of which contain edges labeled $\{1, 2\}$, it follows that all the 8 edges labeled $\{2, 3\}$ or $\{4, 1\}$ are in the remaining family. However, this is impossible because by Lemma 2.3(3) each family contains at most $n_2/2 + 2 = 4$ edges. □

**Lemma 6.5.** *Suppose $n_1 \geq 3$, $n_2 \geq 4$, $\Delta \geq 4$, and suppose $\Gamma_1^+$ has a full vertex $u_k$. Then*

  *(i) $\widehat{\Gamma}_1^+$ has no vertices of valency at most 1;*

  *(ii) if $u_i$ is a boundary vertex of $\widehat{\Gamma}_1^+$ with valency at most 3 in $\widehat{\Gamma}_1^+$, then $i$ is a label of a Scharlemann cycle in $\Gamma_2$; and*

(iii) *each component of $\widehat{\Gamma}_1^+$ contains at most one boundary vertex of valency at most 3.*

*Proof.* Since $u_k$ is a full vertex, all the $k$-edges in $\Gamma_2$ are negative edges. By Lemma 2.5(4) no three such are parallel on $\Gamma_2$. Since there are $\Delta n_2$ $k$-edges, we see that the reduced graph $\widehat{\Gamma}_2$ contains at least $\Delta n_2/2 \geq 2n_2$ negative edges.

(i) By Lemma 2.12(2) there are at most $3n_2$ edges in the reduced graph $\widehat{\Gamma}_2$. We have shown that at least $2n_2$ of the edges are negative, so $\Gamma_2$ has at most $n_2$ positive edges. By Lemma 2.3(3) each edge of $\widehat{\Gamma}_1^+$ represents at most $n_2/2 + 2 \leq n_2$ edges of $\Gamma_1$, so if some $u_i$ has valency at most 1 in $\widehat{\Gamma}_1^+$, then it has more than $2n_2$ negative edges, so there are more than $2n_2$ positive $i$-edges in $\Gamma_2$, which then implies that 3 of those edges are parallel in $\Gamma_2$. This contradicts Lemma 2.5(4).

(ii) If $u_i$ has more than $n_2$ negative edges, then two of them are parallel positive edges on $\Gamma_2$, hence by Lemma 2.4 $i$ is a label of some Scharlemann cycle in $\Gamma_2$, and we are done. Now assume that $u_i$ has at most $n_2$ negative edges, and hence at least $(\Delta - 1)n_2$ positive edges. By assumption it has at most three families of parallel positive edges, and by Lemma 2.3(3) each family contains at most $n_2/2 + 2$ edges, so we have $(\Delta - 1)n_2 \leq 3(n_2/2 + 2)$, which implies that $\Delta = n_2 = 4$, and each family contains exactly $n_2/2 + 2 = 4$ edges. Since $u_i$ is a boundary vertex in $\widehat{\Gamma}_1^+$, these three families are adjacent, as shown in Figure 6.1.

Figure 6.1

By Lemma 2.2(6) $\Gamma_1$ contains no extended Scharlemann cycle, so the labels at the other endpoints of the positive edges are as shown in Figure 6.1. From the figure we can see that there are three Scharlemann cycles with label pair $\{1, 2\}$ and another three with label pair $\{3, 4\}$. On $\Gamma_2$, this implies that there are 6 $i$-edges connecting $v_1$ to $v_2$, and another 6 connecting $v_3$ to $v_4$. Since $n_2 = 4$, there are no other vertices, so it is easy to see that there exist three parallel $i$-edges, contradicting Lemma 2.5(4).

(iii) If $G$ is a component of $\widehat{\Gamma}_1^+$, and $u_i, u_j$ are boundary vertices of valency at most 3 in $G$, then by (ii) both $i, j$ are labels of Scharlemann cycles in $\Gamma_2$. Since $u_i, u_j$ are parallel, no label pair of a Scharlemann cycle could contain both $i$ and $j$; hence there must be two Scharlemann cycles in $\Gamma_2$ with distinct label pairs, contradicting Lemma 2.5(1). □

Consider a graph $\Gamma$ on a sphere $S$. A component $G$ of $\Gamma$ is called an *extremal component* if there is a disk $D$ of $S$, such that $D \cap \Gamma = G$.

**Lemma 6.6.** *Suppose $n_1 \geq 3$, $n_2 \geq 4$, and $\Delta \geq 4$. Then $\Gamma_1$ has no full vertices.*

*Proof.* Suppose $\Gamma_1$ has a full vertex $u_k$. Embed the annulus $A$ in a sphere $S$, and consider $\widehat{\Gamma}_1$ as a graph on $S$. If $\widehat{\Gamma}_1^+$ is connected, then all vertices of $\Gamma_1$ are full. By Lemma 2.12(1) there is a vertex of valency at most 5 in $\widehat{\Gamma}_1^+$, which would contradict Lemma 6.4. Therefore $\widehat{\Gamma}_1^+$ is not connected, so it has at least two extremal components $G_1, G_2$. Let $D_i$ be a disk on $S$ such that $D_i \cap \widehat{\Gamma}_1^+ = G_i$. By Lemma 6.5(1) all vertices of $G_i$ have valency at least 2. That is, $G_i$ is not an isolated vertex, and $a_1(G_i) = 0$.

If some $D_i$ contains no component of $S - A$, then $G_i$ has no loops or parallel edges, so by Lemma 6.2 it has at least two vertices of valency at most 3, which contradicts Lemma 6.5(iii). It follows that each $D_i$ contains exactly one component of $S - A$, hence $l(G_i) + p(G_i) \leq 1$. This also shows that $G_1, G_2$ are the only extremal components, otherwise choosing another extremal component as $G_1$ would give a contradiction. Applying Lemma 6.2 to $G_i$, we have $\sigma(G_i) = 3a_1 + 2a_2 + a_3 + 4l + 2p \geq 6$. Since $a_1 = 0$, and $l + p \leq 1$, we have either $a_3 \geq 2$ or $a_2 \geq 1$. The first case is impossible by Lemma 6.5(iii). Hence $a_2 \geq 1$, that is, each $G_i$ contains a boundary vertex $u^i$ of valency 2. By Lemma 6.5(ii), the index of $u^i$ is a label of all Scharlemann cycles in $\Gamma_2$. Without loss of generality, we may assume that the vertex $u^i$ has label $i$, that is, $u^i = u_i$, $i = 1, 2$. Thus all the Scharlemann cycles in $\Gamma_2$ have label pair $\{1, 2\}$.

CLAIM. *There are exactly 4 Scharlemann cycles of length 2 in $\Gamma_2$, and there are no other positive edges in $\Gamma_2$.*

Let $k$ be the number of positive $i$-edges in $\Gamma_2$, $i = 1, 2$. These are the negative edges in $\Gamma_1$ incident to $u_i$. Each edge of $\widehat{\Gamma}_1^+$ incident to $u_i$ represents at most $n_2/2 + 2$ parallel edges. Since $u_i$ has valency 2 in $\widehat{\Gamma}_1^+$, we have $k \geq \Delta n_2 - 2(n_2/2 + 2) \geq 3n_2 - 4$. On the other hand, as shown in the proof of Lemma 6.5(i), the existence of a full vertex of $\widehat{\Gamma}_1^+$ implies that there are at most $n_2$ positive edges in the reduced graph $\widehat{\Gamma}_2$, and by Lemma 2.5(4) each such edge represents at most two $i$-edges. Therefore, we must have $k \leq 2n_2$. These inequalities give $n_2 = \Delta = 4$, and $k = 2n_2 = 8$. Moreover, each of the 4 positive edge in $\widehat{\Gamma}_2$ must contain exactly two $i$-edges, hence the 8 $i$-edges form 4 pairs of Scharlemann cycles. As shown before the Claim, they all have the label pair $\{1, 2\}$.

By Lemma 2.2(5) the edges of a Scharlemann cycle form an essential loop in $\Gamma_1$, which separates the two extremal components $G_1, G_2$ because $u_1, u_2$ are boundary vertices and neither $G_i$ is contained in a disk in the annulus $A$. Since all the 8 negative edges from $u_1$ end up on $u_2$, there are no other negative edges in $\Gamma_1$. This completes the proof of the claim. In particular, all vertices of $\Gamma_1$ other than $u_1, u_2$ are full vertices.

Note that since there is an arc connecting the only two extremal components $G_1, G_2$ of $\widehat{\Gamma}_1^+$, we must have $\widehat{\Gamma}_1^+ = G_1 \cup G_2$. We have assumed that $n_1 \geq 3$, hence one of the $G_i$, say $G_1$, contains more than one vertex.

Embed $A$ in a sphere $S$ and consider $G_1$ as a graph on $S$. Let $V, E, F$ be the number of vertices, edges and faces of $G_1$. By Lemma 6.4 all the vertices but $u_1$ have valency at least 6, so
$$2 + 6(V - 1) \leq 2E.$$
There are two faces $D_1, D_2$ of $G_1$ which may contain some components of $S - A$, all others are incident to at least three edges. Since $u_1$ has valency 2, and $G_1$ is connected and contains at least two vertices, there is no loop based at $u_1$, so the disk face $D_i$ containing $u_1$ on its boundary has at least two edges. Hence we have
$$2 + 1 + 3(F - 2) \leq 2E.$$
Solving these two inequalities, we have $3V - 3E + 3F \leq 5$, which contradicts the Euler number formula $V - E + F \leq 2$. □

**Lemma 6.7.** *Suppose $n_1 \geq 3$, $n_2 \geq 4$, and $\Delta \geq 4$. If $\widehat{\Gamma}_1^+$ has no interior vertex, then $n_2 = \Delta = 4$. Moreover, each vertex of $\Gamma_1$ has at most 8 negative edges incident to it, and there are 8 negative edges connecting a pair of vertices $u_1, u_2$, which form 4 Scharlemann cycles in $\Gamma_2$.*

*Proof.* Let $u_i$ be a vertex which has valency 2 in $\widehat{\Gamma}_1^+$, given by Lemma 6.2(i). By Lemma 2.3(3) at $u_i$ there are at most $2(n_2/2 + 2) = n_2 + 4$ endpoints of positive edges, hence at least $(\Delta - 1)n_2 - 4$ negative edges. Choose a vertex $u_j$ such that the number $m$ of negative edges incident to it is maximal among all the vertices in $\Gamma_1$. The above shows that $m \geq (\Delta - 1)n_2 - 4$. The parity rule says that there are $m$ positive $j$-edges in $\Gamma_2$.

Let $G$ be the subgraph of $\Gamma_2$ consisting of the positive $j$-edges and all vertices of $\Gamma_2$. Then $G$ has $n_2$ vertices and $m$ edges. Denote by $F$ the number of disk faces of $G$. Calculating the Euler number, we have $n_2 - m + F \geq 0$, hence
$$F \geq m - n_2.$$
By Lemma 2.11, each disk face of $G$ contains a Scharlemann cycle, so there are at least $t = 2F \geq 2(m - n_2)$ edges of Scharlemann cycles. By Lemma 2.5(1), all these Scharlemann cycles have the same label pair, say $\{1, 2\}$. Let $p \in \{1, 2\}$. Then in $\Gamma_1$ there are at least $t$ negative edges incident to the vertex $u_p$. By the definition of $m$, $u_p$ also has at least $\Delta n_2 - m$ endpoints of positive edges. Thus the valency of $u_p$ is at least
$$\begin{aligned}
t + (\Delta n_2 - m) &\geq 2(m - n_2) + \Delta n_2 - m = m - 2n_2 + \Delta n_2 \\
&\geq (\Delta - 1)n_2 - 4 - 2n_2 + \Delta n_2 \\
&= (\Delta - 3)n_2 - 4 + \Delta n_2 \\
&\geq \Delta n_2
\end{aligned}$$

Since the valency of $u_p$ is $\Delta n_2$, all the above inequalities are equalities, so we must have (i) $\Delta = n_2 = 4$, (ii) $m = (\Delta - 1)n_2 - 4 = 8$, hence there are at most 8 negative edges at each vertex in $\Gamma_2$, and (iii) there are $t = 2(m - n_2) = 8$ edges of Scharlemann cycles with label pair $\{1, 2\}$, hence on $\Gamma_1$ there are exactly 8 negative edges at $u_i$, $i = 1, 2$, all connecting $u_1$ to $u_2$. $\square$

**Lemma 6.8.** *Suppose $n_1 \geq 3$, $n_2 \geq 4$, $\Delta \geq 4$, and suppose $\widehat{\Gamma}_1^+$ has no interior vertex. If $\Gamma_1$ has at least three vertices incident to 8 negative edges, then $\Gamma_1$ is 2-separable.*

*Proof.* We have seen from Lemma 6.7 that $\Delta = n_2 = 4$, and $u_1, u_2$ have 8 negative edges, which form 4 Scharlemann cycles in $\Gamma_2$ with label pair $\{1, 2\}$. Assume $u_i$ is another vertex incident to 8 negative edges. By Lemma 2.5(1) $i$ is not a label of a Scharlemann cycle, so by Lemma 2.4 no two of the 8 positive $i$-edges in $\Gamma_2$ are parallel. Thus the reduced graph $\widehat{\Gamma}_2$ has at least 8 positive edges. By Lemma 2.12 $\widehat{\Gamma}_2$ has at most $3n_2 = 12$ edges, so it has at most $12 - 8 = 4$ negative edges.

For each label $j$, by Lemma 2.5(4) each of the 4 negative edges in $\widehat{\Gamma}_2$ contains at most two $j$-edges, so $\Gamma_2$ has at most 8 negative $j$-edges. On $\Gamma_1$ this implies that there are at most 8 positive edges incident to $u_j$. On the other hand, by Lemma 2.7 there are also at most 8 negative edges at $u_j$. Since there are a total of 16 edges at each vertex, this implies that there are exactly 8 positive and 8 negative edges at each vertex $u_j$. In particular, both $\Gamma_\alpha$ have the same number of positive edges and negative edges. There are a total of $\Delta n_1 n_2/2 = 8n_1$ edges, so $\Gamma_2$ has $4n_1$ negative edges. By Lemma 2.5(3) each of the 4 negative edges in $\widehat{\Gamma}_2$ represents exactly $n_1$ negative edges.

Let $\mathcal{E}$ be a set of $n_1$ parallel negative edges of $\Gamma_2$. By Lemma 2.2(3) they form several essential cycles of length $L$ on $\Gamma_1$, which we will call the $\mathcal{E}$-cycles. Choose $\mathcal{E}$ so that $L$ is maximal. Let $e_1, e_2$ be a Scharlemann cycle in $\Gamma_2$, with label pair $\{1, 2\}$. By Lemma 2.2(5) these edges form an essential cycle in $\Gamma_1$, hence cut the annulus $A$ into two annuli $A_1, A_2$, each containing several $\mathcal{E}$-cycles, and having $u_1, u_2$ on its boundary. Let $p$ (resp. $n$) be the number of $\mathcal{E}$-cycles in $A_1$ whose vertices are positive (resp. negative). Assume that $u_1$ is a positive vertex. Then on $A_1$ there are $pL$ positive vertices, and $nL + 1$ negative vertices (including $u_2$). Recall that each vertex has 8 negative edges, which connect it to vertices of opposite sign. Moreover, all the negative edges at $u_1$ go to $u_2$. Counting the number of negative edges, we have

$$8(pL) = 8(nL + 1),$$

which implies that $L = 1$.

When $L = 1$, all the positive edges in $\Gamma_1$ are loops. Thus the loops at $u_1, u_2$ and the edges connecting them form a component of $\Gamma_1$, contained in an annulus $A'$ with boundary disjoint from $\Gamma_1$. Therefore $\Gamma_1$ is 2-separable. $\square$

**Lemma 6.9.** *Suppose $n_1 \geq 3$, $n_2 \geq 4$, and $\Delta \geq 4$. If $\widehat{\Gamma}_1^+$ has no full vertices, then it is 2-separable.*

*Proof.* By Lemma 6.7, $n_2 = \Delta = 4$, and there are four Scharlemann cycles on $\Gamma_2$, all having the same label pair $\{1, 2\}$, say. Let $G_i$, $i = 1, 2$, be the component of $\widehat{\Gamma}_1^+$ containing $u_i$. If there is an extremal component $G$ of $\widehat{\Gamma}_1$, $G \neq G_1, G_2$, or if $G_i$ is an extremal component but $u_i$ is a cut vertex, then applying Lemma 6.2(ii) to $G_1, G_2$ and $G$ we see that there is a vertex $u_j \neq u_1, u_2$, which has valency at most 2 in $\widehat{\Gamma}_1^+$. Such a vertex has at most $2(n_2/2+2) = 8$ positive edges, hence at least $\Delta n_2 - 2(n_2/2+2) = 8$ negative edges, so the result will follow from Lemma 6.8. Therefore we assume that $G_1, G_2$ are the only extremal components of $\widehat{\Gamma}_1^+$, and $u_1, u_2$ are boundary vertices. Since there are edges connecting $u_1$ to $u_2$, there are no other components, that is, $\widehat{\Gamma}_1^+ = G_1 \cup G_2$.

The edges $e_1, e_2$ of a Scharlemann cycle in $\Gamma_2$ form an essential cycle on the annulus $A$, which cuts $A$ into two annuli $A_1, A_2$. Since $u_i$ is not a cut vertex, each $G_i$ lies in one of the $A_j$. The two $G_i$ must be on different $A_j$, otherwise one of the $G_i$ would lie on a disk in $A$, so by Lemma 6.2(ii) it would have at least two vertices of valency at most 2 in $\widehat{\Gamma}_1^+$. Together with the vertices $u_1$ or $u_2$ on the other component $G_k$, this would give three vertices incident to at least 8 negative edges, contradicting Lemma 6.8.

Now the only edges connecting $G_1$ to $G_2$ must have at least one endpoint on $u_1$ or $u_2$; but since all the 8 negative edges from $u_1$ must go to $u_2$ (Lemma 6.7), there are no other negative edges. It follows that all vertices other than $u_1, u_2$ are full vertices. Since $n_1 \geq 3$, this contradicts the assumption.  $\square$

**Theorem 6.10.** *Suppose $M$ is a hyperbolic manifold, such that $M(r_1)$ contains an essential annulus $\widehat{F}_1$, $M(r_2)$ contains an essential torus $\widehat{F}_2$, and $\Delta(r_1, r_2) \geq 4$.*

*(1) If $\Delta = 4$, $M$ is the exterior of the link in Figure 7.1(a) or (b); otherwise, $\Delta = 5$, and $M$ is the exterior of the link in Figure 7.1(c).*

*(2) $F_\alpha$ can be chosen to intersect the Dehn filling solid torus at most twice.*

*(3) If some $F_\alpha$ intersects the Dehn filling solid torus only once, then $M$ is the exterior of the Whitehead link in Figure 7.1(a).*

*Proof.* If some $n_\alpha = 1$, then by Proposition 3.5 the manifold $M$ is the exterior of the Whitehead link in Figure 7.1(a), $n_\beta = 2$, and $\Delta = 4$.

Assume $n_1, n_2 \geq 2$, and some $n_\alpha = 2$. If $n_1 = 2$, then by definition $\Gamma_1$ is 2-separable. If $n_2 = 2$ and $n_1 > 2$, then by Proposition 4.6 $\Gamma_1$ is also 2-separable. By Proposition 5.5, in this case we must have $n_1 = n_2 = 2$, and the manifold $M$ is the exterior of the link in Figure 7.1(b) or (c) respectively, according to $\Delta = 4$ or 5.

The remaining case is that $n_\alpha \geq 3$ for $\alpha = 1, 2$. By Lemma 6.1, we cannot have $n_2 = 3$. Assume $n_2 \geq 4$. By Lemma 6.6 the graph $\widehat{\Gamma}_1^+$ has no full vertex, and by Lemma 6.9 this implies that $\Gamma_1$ is 2-separable, which by Proposition 5.5 implies that $n_1 = 2$, contradicting the assumption. Hence this case does not happen.  $\square$

## §7. Some Examples

In this section we give some examples of manifolds admitting toroidal and annular surgeries. Consider the three links shown in Figure 7.1. The first one is the Whitehead link. It is the 2-bridge link associated to the rational number 3/8. The second is the 2-bridge link associated to 3/10. The third one is called the Whitehead sister link.

Figure 7.1

It is the $(-2, 3, 8)$ pretzel link. The Whitehead link is well known to be hyperbolic; see for example [Th2]. We will show in Theorem 7.5 that the other two links are also hyperbolic. First consider the Whitehead link in Figure 7.1(a).

**Lemma 7.1.** *Let $M$ be the exterior of the Whitehead link $L$ shown in Figure 7.1(a). Then there are two slopes $r_1, r_2$ such that (i) $M(r_1)$ contains an essential torus and an essential annulus, each intersecting the Dehn filling solid torus $J_1$ once; (ii) $M(r_2)$ contains an essential annulus and an essential torus, each intersecting $J_2$ twice; and (iii) $\Delta(r_1, r_2) = 4$.*

*Proof.* The manifold $M$ is shown in Figure 7.2(a), with $T_0 \subset \partial M$ the inside torus. Let $\varphi$ be the $\pi$-rotation of $(M, T_0)$ about the horizontal axis shown in Figure 7.2(a). Let $\alpha$ be the axis of $\varphi$ in $M$. Then one can see that $(M, \alpha)/\varphi = (X, \eta)$, where $X$ is a product $S^2 \times I$, and $\eta$ consists of four arcs. See Figure 7.2(b). Thus $M$ is the double cover of $X$ branched along $\eta$. By an isotopy $(X, \eta)$ is equivalent to the pair shown in Figure 7.2(c).

It is well known that, with respect to a certain framing on $T_0$, the manifold $M(r)$ can be obtained by filling the inner component of $X$ with a rational tangle to get a pair $(B, \eta(r))$, and then taking the double cover of $B$ branched along $\eta(r)$; see [M]. The tangle $(B, \eta(1/0))$ is shown in Figure 7.2(d). Consider the shaded disk $D$ and Möbius band $U$ shown in Figure 7.2(d). Cutting along $D$ and $U$, we get a manifold $X$ homeomorphic to $T \times I$, on which $D$ becomes an annulus $A$, and $U$ becomes a torus $T$.

The double branched cover $M(-1/2)$ can be obtained by taking two copies of $X$, gluing $T$ to $T$, and $A$ to $A$. Since $T$ and $A$ are incompressible in both sides, they become essential surfaces in $M(-1/2)$. The core of the Dehn filling solid torus $J$ branch covers the core of the rational tangle, shown as $K$ in Figure 7.2(d). Since each of $D$ and $U$ intersects $K$ in a single boundary point, the annulus $A$ and the torus $T$ also intersect $J$ in a single disk. This proves (i).

Now consider the tangle $(B, \eta(1/2))$ shown in Figure 7.2(e). It is equivalent to that in Figure 7.2(f). There is a disk $D$ cutting it into a $(-1/2, -1/2)$-Montesinos tangle and a $(1/2)$-rational tangle. The lift of $D$ is then an essential annulus $A$ in $M(1/2)$, cutting $M(1/2)$ into a Seifert fiber space $Y$ whose orbifold is a disk with two cone points of order 2, and a solid torus on which $A$ runs twice along the longitude. Thus $A$ is an essential annulus, and the boundary of $Y$ pushed into the interior of $M(1/2)$ is an essential torus, which covers a Conway sphere $S$ on the right hand side of the disk $D$ in $B$. Notice that both $D$ and $S$ intersect the core of the rational tangle just once, hence $A$ and $T$ intersect the core of the Dehn filling solid torus exactly twice. Since $\Delta(-1/2, 1/2) = 4$, this completes the proof of (ii) and (iii).

Notice that the framing we have chosen is not the standard meridian-longitude framing of the Whitehead link. By a more careful examination, one can show that $M(-1/2)$ is the 0-surgery, and $M(1/2)$ is the 4-surgery with respect to the standard framing. $\square$

Figure 7.2

**Lemma 7.2.** *Suppose $M$ is a hyperbolic manifold such that $M(r_1)$ is annular, $M(r_2)$ is toroidal, and the intersection graphs $\Gamma_\alpha$ are as shown in Figure 3.1. Then $M$ is the exterior of the Whitehead link.*

*Proof.* One can prove this by using the argument in the proof of Lemma 3.2. Here we give an explicit construction of $M$. The method can be used to construct other manifolds from intersection graphs.

Let $J_1$ be the Dehn filling solid torus in $M(r_1)$. Taking the union of $A \times I$ with $J_1$, so that the fat vertex of $\Gamma_1$ in $A$ is identified with a meridian disk of $J_1$, we get a genus 2 handlebody $X$, as shown in Figure 7.3(a). The graph $\Gamma_2$ in Figure 3.1(b) cuts the torus $T$ into two disk faces, whose boundary curves are parallel on $\partial X$, each being isotopic to the curve $C$ shown in Figure 7.3(a). Therefore the regular neighborhood of the union of the handlebody $X$ and the punctured torus $F_2$, with the spherical boundary component capped off by a 3-ball, is the same as the manifold $Y$ obtained by attaching a 2-handle $H$ to $X$ along the curve $C$. Let $K_1$ be the central curve of the solid torus $J_1$. The manifold $M$ is obtained by removing a regular neighborhood $V$ of $K_1$ from $Y$.

There is an involution $\varphi$ of $X$, as shown in the figure, which carries $C$ to itself, so $\varphi$ extends to an involution of $Y = X \cup H$. Let $L$ be the axis of $\varphi$ in $Y$. After the involution, the quotient of $X$ is a 3-ball $B = X/\varphi$, and the image of $L \cap X$ consists of three arcs $\alpha_1, \alpha_2, \alpha_3$ in $B$, as shown in Figure 7.3(b). The quotient of $C$ is an arc $\alpha$ on $\partial B$ connecting the two endpoints of $\alpha_3$. The image of the attached 2-handle $H$ is a 3-ball $B'$, with $L \cap H$ projecting to an arc which is isotopic rel $\partial$ to $\alpha$. It follows that there is a homeomorphism from $Y/\varphi$ to $B$, under which the image of $L$ is the union of $\alpha$ and the three arcs $\alpha_i$, with $\alpha$ pushed into the interior of $B$. Denote by $\beta$ the 1-manifold $L/\varphi$. The central curve $K_1$ of $J_1$ becomes an arc with endpoints on $\beta$, denoted by $K$. In Figure 7.3(b) $K$ is the thick grey arc inside $B$. Consider $\beta$ as a tangle in the 3-ball $B$. After an isotopy, the triple $(B, \beta, K)$ is the same as that shown in Figure 7.2(b), where $B$ is the 3-ball, and $\beta$ is the tangle consisting of two arcs and a circle.

Let $X = B - \text{Int}N(K)$, and $\eta = \beta \cap X$. Then $M$ is the double cover of $X$ branched along $\eta$, which has been shown in the proof of Lemma 7.2 to be the exterior of the Whitehead link. □

Figure 7.3

We need a lemma to test the hyperbolicity of the exteriors of braids in solid tori.

A braid is *trivial* if it has only one strand.

**Lemma 7.3.** *Let $K$ be a knot in a solid torus $V$, which is a nontrivial braid, and let $M$ be the exterior of $K$ in $V$. Then $M$ is hyperbolic unless (i) there is a solid torus $V'$ in $V$ whose core is a nontrivial braid in $V$, and $K$ is a nontrivial braid in $V'$; or (ii) $K$ is a torus knot in $V$, i.e. it is isotopic to a curve on $\partial V$.*

*Proof.* Clearly $M$ is irreducible and $\partial$-irreducible. Thus by the Geometrization Theorem of Thurston for Haken manifolds [Th1], $M$ is hyperbolic unless it is toroidal or Seifert fibered.

Suppose $M$ is toroidal, containing an essential torus $T$. Let $D$ be a meridian disk of $V$ intersecting $K$ in $n$ disks, where $n$ is the braid index of $K$. Let $P$ be the punctured disk $D - \text{Int} N(K)$. Isotope $T$ so that it has minimal intersection with $P$. After cutting along $P$, $M$ becomes a product $P \times I$, and $T$ becomes a set of essential annuli $\mathcal{A}$ with boundary on $P \times \partial I$. Hence $\mathcal{A}$ is a product $\mathcal{C} \times I$, where $\mathcal{C}$ is a set of essential curves on $P$. It is now clear that after gluing $P \times 0$ to $P \times 1$, the annuli $\mathcal{A}$ are glued to give a torus bounding a solid torus $V'$ whose core is a braid in $V$, which is nontrivial because $T$ is not boundary parallel. $K$ is now a braid in $V'$, which must be nontrivial because $T$ is not parallel to $\partial N(K)$. Hence (i) holds in this case.

Now suppose $M$ is Seifert fibered. After filling in $N(K)$, the Seifert fibration extends to one on the solid torus $V$, whose orbifold is a disk with at most one cone point. Since $K$ is nontrivial, it is not a singular fiber, hence it is isotopic to a regular fiber on $\partial V$, and (ii) follows. □

A braid $b$ can be written as $\sigma_{i_1}^{\epsilon_1} \sigma_{i_2}^{\epsilon_2} \ldots \sigma_{i_k}^{\epsilon_k}$, where each $\sigma_{i_j}$ is one of the standard generators of the braid group, and $\epsilon_j = \pm 1$. The *total exponent* of $b$ is defined as $e(b) = \sum \epsilon_j$.

**Corollary 7.4.** *Suppose a knot $K$ is an $n$-braid in a solid torus $V$, with total exponent $e$. If $n$ is prime, and if $e$ is not a nonzero multiple of $n-1$, then the exterior of $K$ is hyperbolic.*

*Proof.* If $K$ is an $n_1$ braid in $V'$ whose core is an $n_2$ braid in $V$, then $K$ is an $n_1 n_2$ braid in $V$. If $K$ is a $(n, q)$ torus knot in $V$, then the total exponent is $(n-1)q$. Hence the corollary follows from Lemma 7.3. □

Now let $M_i$, $i = 1, 2, 3$ be the exterior of the three links $L_i$ shown in Figure 7.1. Denote by $T_0$ the boundary torus corresponding to the right component of $L_i$ in Figure 7.1.

**Theorem 7.5.** *Each $M_i$ is a hyperbolic manifold, and admits two Dehn fillings $M_i(r_1)$ and $M(r_2)$, such that*

  (i) $\Delta = 4$ *for* $i = 1, 2$, *and* $\Delta = 5$ *for* $i = 3$;
  (ii) $M_1(r_1)$ *contains an essential annulus $A$ and an essential torus $T$, each intersecting $J_1$ once. All the other $M_i(r_k)$ contain an essential annulus $A$ and an essential torus $T$, each intersecting $J_k$ twice.*

*Proof.* If $i = 1$, this is Lemma 7.1. So assume $i = 2$ or 3.

Removing the left component of $L_i$, we get a solid torus $V$ containing a knot $K_i$. For $K_2$, the braid index is 3, and the total exponent is 0. For $K_3$, the braid index is 5, and the total exponent is 6. It follows from Corollary 7.4 that the exteriors of both $K_i$ are hyperbolic.

Figure 7.4

The proof of the rest part of the theorem is similar to that of Lemma 7.1, with Figure 7.2 replaced by Figures 7.4 and 7.5. The essential torus and essential annulus in $M_2(r_2)$ are the double covers of the disk $D_1$ in Figure 7.4(d) and the boundary sphere $S_1$ on the right hand side of $D$, pushed slightly into the interior. Similarly for the others. In all cases, the tangle $(B, \eta(r_j))$ is the sum of a $(\pm 1/2, \mp 1/3)$ Montesinos tangle and a $(1/2)$-rational tangle, so $M_i(r_j)$ is a graph manifold obtained by gluing a solid torus to a Seifert fiber space along an annulus which runs twice along the longitude of the solid torus. The disks $D_i$ and the Conway spheres $S_i$ intersect the core of the rational tangle at one interior point, so after double covering the annuli and the tori intersect the Dehn filling solid torus at two points. $\square$

Figure 7.5

Consider the manifold $M$ shown in Figure 7.6(a). It is the exterior of a three component link in $S^3$. Let $L = l_0 \cup l_1$ be the link in the solid torus $V$, where $l_0$ is the component with winding number 1 in $V$. Let $T_0 = \partial N(l_0)$, with the standard meridian-longitude framing.

**Lemma 7.6.** *$M$ is hyperbolic, and both $M(0)$ and $M(-3/2)$ are toroidal and annular.*

*Proof.* $L$ is a braid in $V$, so $M$ is irreducible and $\partial$-irreducible. As in the proof of Lemma 7.3, it can be shown that if $T$ is an essential torus in $M$, then $T$ bounds a solid torus $V'$ in $V$ such that the core of $V'$ is a braid in $V$, and $L \cap V'$ is a braid in $V'$. If $V'$ contains the component $l_0$ which has braid index 1 in $V$, then the core of $V'$ is a trivial braid in $V$. So if $V'$ also contains $l_1$ then $T = \partial V'$ is parallel to $\partial V$, and if $V'$ does not contain $l_1$ then $T$ is parallel to $\partial N(l_0)$, in either case contradicting the assumption that $T$ is essential in $M$. Hence $V'$ must contain $l_1$ but not $l_0$. The braid index of $l_1$ in $V'$ is 2, otherwise $T$ would be parallel to $\partial N(l_1)$. Let $K'$ be the core of $V'$. Then the linking number $lk(l_0, l_1) = 2\, lk(l_0, K')$ is even, which is absurd because from the picture we see that $lk(l_0, l_1) = 1$.

To prove that $M$ is hyperbolic, it remains to show that $M$ is not a Seifert fiber space. If $M$ were Seifert fibered, then after the trivial filling on $L$, the Seifert fibration would extend to one on the solid torus $V$, so that each $l_i$ is a fiber. Since a Seifert fibration of $V$ has at most one singular fiber, and $l_0, l_1$ are not isotopic, one can see that $l_0$ must be a singular fiber while $l_1$ is a regular fiber, so there would be an isotopy sending $l_0$ to the core $l'_0$ of $V$ and $l_1$ to a $(2,1)$ curve $l'_1$ on $\partial V$. However, the braid index of $l'_0 \cup l'_1$ is 3 but the braid index of $l_0 \cup l_1$ is 1, so they cannot be isotopic. This completes the proof that $M$ is hyperbolic.

As before, let $\varphi$ be the $\pi$-rotation of $M$ about the horizontal axis shown in Figure 7.6(a), and let $\alpha$ be the axis of $\varphi$ in $M$. Then $(M, \alpha)/\varphi = (X, \eta)$, where $X$ is a twice punctured 3-ball, and $\eta$ consists of 6 arcs, as shown in Figure 7.6(b). $M$ is the double cover of $X$ branched along $\eta$. After an isotopy, the picture is shown in Figure 7.6(c), in which the image of $T_0$ is the left sphere $S_0$. We have been careful to ensure that the standard meridian and longitude framings induce the $\infty$ and $0$ slopes on $\partial S_0$. However, notice that the conventions about slopes on $\partial S_0$ and $T_0$ are opposite. Denote by $(Y, \eta(r))$ the pair obtained from $(X, \eta)$ by filling $S_0$ with a rational tangle of slope $r$. Then $M(r)$ is the double cover of $(Y, \eta(-r))$.

Figure 7.6

It is easy to see that $(Y, \eta(0))$ is the tangle shown in Figure 7.6(d). A disk intersecting $\eta(0)$ twice is shown in the figure. The double branched cover $M(0)$ of $(Y, \eta(0))$ is the union of a trefoil knot exterior and a copy of $T^2 \times I$, glued along an annulus which is meridional on the trefoil knot exterior. Hence $M(0)$ is actually the exterior of the connected sum of a trefoil and a Hopf link. This manifold is toroidal and annular.

Figure 7.6(e) is the picture of $(Y, \eta(3/2))$. After an isotopy it is shown in Figure 7.6(f). In this case $M(-3/2)$ is the union of a solid torus and a $(2,1)$ cable space, glued along an annulus which runs twice along the longitude of the solid torus, hence it is toroidal and annular. One can show that $M(-3/2)$ is actually the exterior of a knot $K$ in a solid torus $V$, where $K$ is some $(2,q)$ cable of a $(2,1)$ torus knot in $V$. □

We call the manifold $M$ above the *magic manifold*. Table 7.7 lists some manifolds obtained by Dehn filling on $T_0$. We use $X(p/q)$ to denote the exterior of the 2-bridge link associated to the rational number $p/q$, and $X(p,q,r)$ the exterior of the $(p,q,r)$ pretzel link. We have seen from Lemma 7.6 that $M(0)$ and $M(-3/2)$ are toroidal and annular. The first is the exterior of a link which is the connected sum of a trefoil and a Hopf link, and the second one is the exterior of some $(2,q)$ cable of the $(2,1)$ torus knot in the solid torus, so it is also the exterior of a link in $S^3$.

The first five fillings are non-hyperbolic, and the last three are the exteriors of the three links in Figure 7.1, which by Theorem 6.10 are the only manifolds admitting annular and toroidal fillings of distance at least 4 apart. It is surprising that all the three exceptional manifolds in Theorem 6.10 arise from Dehn fillings on this manifold.

Table 7.7

There is an automorphism of $M$ interchanging any pair of boundary components of $M$. Let $T_1$ be the boundary of the component of $L$ which wraps along $V$ twice. Let $(m_i, l_i)$ be the standard meridian-longitude pair of $T_i$. The homeomorphism sending $T_0$ to $T_1$ will send $m_0$ to $l_1$, and $l_0$ to $m_1 - l_1$. Hence all the above Dehn fillings can also be realized by Dehn filling on $T_1$, and the slopes can be calculated accordingly. For example, the two toroidal and annular fillings on $T_1$ have slopes $-1$ and $2$, respectively.

We can apply Theorem 6.10 to prove the following result.

**Proposition 7.7.** *If $r \neq 0, -1, -3/2, -2, \infty$, then $M(r)$ is hyperbolic.*

*Proof.* Note that $M$ has three boundary components, so it is not one of the three exceptional manifolds in Theorem 6.10. Let $r_1 = 0$, and $r_2 = -3/2$. By Lemma 7.6 each $M(r_i)$ is toroidal and annular, so if $r$ is a slope such that $\Delta(r, r_i) \geq 4$ for some $i$, then by Theorem 6.10 the manifold $M(r)$ is atoroidal and anannular, and by [Wu1,Oh] and [GLu] it is also irreducible and $\partial$-irreducible, and hence $M(r)$ is hyperbolic. Thus if $M(r)$ is non-hyperbolic, then $\Delta(r, r_i) \leq 3$ for $i = 1, 2$.

Suppose $r = a/b$. Without loss of generality we may assume that $a \geq 0$. Then we have
$$\Delta(0, a/b) = a \leq 3;$$
$$\Delta(-3/2, a/b) = |2a + 3b| \leq 3.$$

Solving these inequalities, we have $a/b = 0, \infty, -1, -3/2, -2$, or $-3$. From Table 7.7 we see that $M(-3)$ is the exterior of the Whitehead sister link shown in Figure 7.1(c), which is hyperbolic by Theorem 7.5. □

**Theorem 7.8.** *There are infinitely many hyperbolic manifolds $M_s$, which admit two Dehn fillings $M_s(r_1), M_s(r_2)$, each containing an essential annulus and an essential torus, with $\Delta(r_1, r_2) = 3$.*

*Proof.* Let $M$ be the manifold in Lemma 7.6. We are now considering Dehn fillings on both $T_0$ and $T_1$, so denote by $M(r, *) = M(r)$ the manifold obtained by $r$-filling on $T_0$, by $M(*, s) = M_s$ the manifold obtained by $s$-filling on $T_1$. Let $r_1 = 0$, and $r_2 = -3/2$.

With respect to certain framings on $T_1$, the manifold $M(*, s)$ can be obtained by filling the sphere $S_1$ in Figure 7.6(c) with a rational tangle of slope $s$, then taking the double branched cover. From Figure 7.6(d) and 7.6(f) we can see that if $s \notin \{0\} \cup \{1/n \,|\, n \in \mathbb{Z}\}$, then for $i = 1, 2$, the essential annulus and essential torus in $M(r_i) = M(r_i, *)$ remain essential in $M(r_i)(s) = M_s(r_i)$. By Thurston [Th2], there is a finite set $\mathcal{S}$ such that if $s \notin \mathcal{S}$ then $M(*, s)$ is hyperbolic. Actually using Proposition 7.7 and the symmetry of $M$ we see that $\mathcal{S}$ contains exactly 5 slopes. Thus if $s \notin \mathcal{S} \cup \{0\} \cup \{1/n \,|\, n \in \mathbb{Z}\}$, then $M_s$ is hyperbolic, and $M_s(r_i)$ is toroidal and annular, $i = 1, 2$. Clearly there are infinitely many such $s$. Since $\Delta(r_1, r_2) = \Delta(0, -3/2) = 3$, this proves the theorem.

Note that the framing on $T_1$ above is not the standard meridian-longitude framing. □

Department of Mathematics, University of Texas at Austin, Austin, TX 78712 and MSRI, 1000 Centennial Drive, Berkeley, CA 94720-5070
*E-mail address*: gordon@math.uiowa.edu

Department of Mathematics, University of Iowa, Iowa City, IA 52242; and MSRI, 1000 Centennial Drive, Berkeley, CA 94720-5070
*E-mail address*: wu@math.uiowa.edu